\documentclass[11pt,a4paper]{article}

\usepackage{amsmath}
\usepackage{amsthm}
\usepackage{amssymb}
\usepackage{latexsym}
\usepackage{graphicx,subfig}
\usepackage{mathrsfs}
\usepackage{pinlabel-test}

\allowdisplaybreaks

\newcommand{\A}{\mathcal A}

\newcommand{\dms}{\textup{diam}(\Sigma)}
\newcommand{\e}{\varepsilon}
\newcommand{\g}{\gamma}
\newcommand{\F}{\mathscr F}
\renewcommand{\SS}{\mathcal G}
\renewcommand{\H}{\mathscr H}
\renewcommand{\k}{k}
\newcommand{\kh}{\kappa_H}
\newcommand{\kg}{\kappa_G}

\renewcommand{\L}{\mathcal L}
\newcommand{\N}{\mathbb{N}}
\newcommand{\R}{\mathbb{R}}
\newcommand{\s}{\sigma}
\renewcommand{\S}{\Sigma}
\renewcommand{\t}{{\tau}}

\newcommand{\V}{\textup{Vol\,}}
\newcommand{\vs}{\textup{Vol\,}(\S)}

\newcommand{\GO}{\textup{(G0) }}
\newcommand{\GI}{\textup{(G1) }}

\newcommand{\ds}{\displaystyle}

\newcommand{\val}{\left((0,1);\R^2\right)}

\newcommand{\weakto}{\rightharpoonup}
\newcommand{\weaksto}{\stackrel{*}{\rightharpoonup}}

\def\Xint#1{\mathchoice
   {\XXint\displaystyle\textstyle{#1}}%
   {\XXint\textstyle\scriptstyle{#1}}%
   {\XXint\scriptstyle\scriptscriptstyle{#1}}%
   {\XXint\scriptscriptstyle\scriptscriptstyle{#1}}%
   \!\int}
\def\XXint#1#2#3{{\setbox0=\hbox{$#1{#2#3}{\int}$}
     \vcenter{\hbox{$#2#3$}}\kern-.5\wd0}}
\def\dashint{\Xint-}

\newtheorem{theo}{Theorem}
\newtheorem{lemma}[theo]{Lemma}
\newtheorem{cor}[theo]{Corollary}
\newtheorem{prop}[theo]{Proposition}

\newtheorem{remark}{\mdseries{\itshape{Remark}}}
\newenvironment{rem}{\begin{remark} \upshape}{\end{remark}}
\newtheorem{definition}{\bfseries{\upshape{Definition}}}
\newenvironment{defi}{\begin{definition} \upshape}{\end{definition}}

\usepackage{color}

\begin {document}

\title{Global minimizers for the doubly-constrained Helfrich energy: the axisymmetric case.}

\author{Rustum Choksi\footnote{Department of Mathematics and Statistics,
McGill University,  Montreal, Canada, {rchoksi@math.mcgill.ca}} \qquad 
Marco Veneroni\footnote{Department of Mathematics ``Felice Casorati", University of Pavia, Pavia, Italy,  {marco.veneroni@unipv.it}}}

\date{\today}

\maketitle
\begin{abstract}
Since the pioneering work of Canham and Helfrich, variational formulations involving curvature-dependent functionals, like the classical Willmore functional, have proven useful for shape analysis of biomembranes. 
We address minimizers of the Canham-Helfrich functional defined over closed surfaces enclosing a fixed volume and 
having fixed surface area. By restricting attention to axisymmetric surfaces, we prove the existence of global minimizers.
\end{abstract}
\medskip

\bigskip

\noindent \textbf{Keywords:} Willmore functional, Helfrich functional, biomembranes, direct method in the calculus of variations.
\medskip

\noindent \textbf{AMS subject classification}: 49Q10, 49J45 (58E99, 53C80).


\section{Introduction and main result}

For compact surfaces $\S$ embedded in $\R^3$, the Canham-Helfrich functional 
is defined by
\begin{equation}
\label{eq:helfrich}
	\H (\S)=  \int_\S \left\{  \frac{\kh}{2}(H-H_0)^2+\kg\, K\right\}dA,
\end{equation}
where the integration is with respect to the ordinary 2-dimensional area measure, $H$ is the sum of the principal curvatures of $\S$, i.e., twice the mean curvature, $K$ is the Gaussian curvature, $\kh,\kg\in \R$ are constant bending rigidities and $H_0\in \R$ is a given spontaneous curvature.

We prove the existence of a global minimizer for $\H$, in the class of finite systems of axisymmetric surfaces, under the constraints that the total area and the total enclosed volume of the surfaces are fixed.
\medskip

Biological membranes and their shapes have attracted  attention from researchers across many areas of mathematics. 
For example,  membranes connect classical problems of differential geometry involving Willmore's functional to studies of shape configurations of biological cells in physics and biology (see, e.g., \cite{Jenkins77}, \cite{JuelicherLipowsky96}). More recently,  researchers in both mathematical analysis and scientific computing have 
directed  efforts to understanding \emph{multiphase} membranes, where phase transitions and pattern formation can be observed (see, e.g.\,\cite{BaumgartHessWebb}, \cite{Fried}, \cite{Du_multi}, \cite{Lowengrub_etal}, \cite{ElliottStinner10}, \cite{Zurlo}).
The modeling of multiphase membranes has  numerous applications associated with  artificial membranes in pharmacology and bioengineering (e.g., \cite{Venkatesan}).

In his seminal work \cite{Canham70}, Canham proposed the functional \eqref{eq:helfrich}, in the case $H_0=0$, in order to model the elastic bending energy of biological membranes formed by a double layer of phospholipids. When immersed in water, these molecules, which are composed by a hydrophilic head and a hydrophobic tail, spontaneously aggregate in order to shield the tails from water, forming a closed bilayer with the heads pointing outwards. Since the thickness of a layer is generally three to four orders of magnitude smaller than the size of the observed cells or vesicles, the bilayer is usually approximated as a two-dimensional surface $\S$ embedded in $\R^3$. 
The functional $\H$ is the most general example of energy which is quadratic in the principal curvatures. The parameter $H_0\in \R$, added by Helfrich \cite{Helfrich73}, accounts for an asymmetry in the composition of the layers and gives rise to a spontaneous curvature of the membrane in absence of other constraints. 
The bending rigidities $\kh$ and $\kg$ are also material-dependent parameters. 
Under the simplifying assumption that the membrane is homogeneous, we choose $H_0,\kh$ and $\kg$  constant. In reality, phases with different levels of aggregation and different rigidities are observed \cite{JuelicherLipowsky93}.  

There are two natural constraints associated with the membrane configuration. 
Since lipid membranes are inextensible, the total area of the membrane should be fixed. On the other hand, 
the membrane is permeable to water but not to dissolved ions.  The resulting osmotic pressure leads then to a constraint on the volume enclosed by the membrane, which can therefore be regarded as constant \cite{JuelicherLipowsky96}. From the point of view of our analysis, the constraint on the area plays a crucial role in obtaining a priori bounds and compactness. The constraint on the volume, instead, does not add any property or difficulty, but it is the combination of these two constraints that makes highly nontrivial the problem of determining the minimizer.

Another important feature of membranes is that they can undergo topological changes, for example, a spherical vesicle can shrink at the equator and eventually split into two vesicles (fission) or a small dome can rise from a point of the surface and grow into a new entity which separates from the original one (budding), see e.g. \cite[Section 3.9]{Seifert97} and \cite{BaumgartDWJ}. In order to be able to describe these kind of phenomena, we do not impose restrictions on the number of components of the minimizers.    
\medskip

Helfrich's functional can be regarded as a generalization of the classical Willmore functional, defined by
$$ W(\S)=\frac 14 \int_\S|H|^2\, dA.$$
In the seminal paper \cite{Simon93}, Leon Simon proved that for each $n\geq 3$ there exists a compact embedded real analytic torus in $\R^n$ which minimizes $W$ among compact embedded surfaces of genus 1. 
Following the direct method of the calculus of variations, he first shows that sequences of minimizers are compact in the sense of measures, and then proves that the limit measure is actually an analytic surface. Simon's proof of regularity relies on the invariance of  Willmore functional under conformal transformations and on the fact that a minimizing surface $\S$ must satisfy $4\pi\leq W(\S)\leq 8\pi$. Owing to the presence of the spontaneous curvature $H_0$ and to the combined area and volume constraints, Helfrich's functional \emph{is not conformally invariant}, and for general values of area and volume, we only know that $0\leq \H$. Therefore, though measure-compactness can be easily transferred to our case, Simon's method for regularity cannot be employed, and we have to find a different approach. Regarding minimization with area and volume constraints in the case $H_0=0$, existence of genus 0  minimizers with fixed isoperimetric ratio was recently proved in \cite{Schygulla}, while \cite{WheelerPR} gives a complete classification of smooth critical points with low energy. 

Existence of minimizers for functionals with weak second fundamental form in $L^2$ was addressed also in \cite{Hutchinson86}, using the theory of varifolds (see also the end of Section \ref{ssec:discussion}). However, in contrast to the mean curvature vector, the scalar mean curvature $H$ does not have a variational characterization and there is no definition of scalar mean curvature for an arbitrary integral varifold. This obstacle can be removed using the  generalized Gauss graphs introduced in \cite{AnzellottiSerapioniTamanini} and developed in \cite{Delladio97}. Compactness and lower-semicontinuity properties allow to obtain a minimizer,  but it is not trivial to understand whether the limit, which in general is only a rectifiable current, is actually a classical surface. This is certainly true in the case of one-dimensional curves in $\R^2$, see for example,  \cite{BellDMasoPaolini93}, \cite{BellettiniMugnai04} and \cite{BellettiniMugnai07}, but the question remains open for surfaces in $\R^3$.

A different approach, based on a new formulation for the Euler-Lagrange equation of Willmore functional, was introduced in \cite{Riviere08}. One of the results therein is a new proof of the existence of minimizers. In the attempt to apply this new method to Helfrich functional, the same difficulties as above appear, in particular, the lack of an equivalent of  Li-Yau minimality condition \cite[Theorem 6]{LiYau82} for Helfrich functional necessitates another approach in order to guarantee that minimizers are embedded. 

Regarding the existence of minimizers for the constrained Willmore functional, in \cite{DallAcqua} it is proven existence and regularity of axisymmetric solutions of the Euler-Lagrange equations with symmetric boundary conditions. In \cite{Schygulla}, adopting the techniques introduced in \cite{Simon93}, the author proves the existence of minimizers with prescribed isoperimetric ratio. We are not aware of any extension of these methods to the constrained Helfrich functional.\medskip

In the present work, we give an answer that is only partial, since we restrict to axisymmetric surfaces. We note that our result cannot be obtained from the above-mentioned results for $W^{2,2}$-regular curves, since if a curve $\g$ generates a surface with bounded Helfrich energy $\H$, it is not true in general that $\g$ is $W^{2,2}$-regular (see Section \ref{ssec:discussion} below). The class of axisymmetric surfaces is probably the most interesting from the point of view of applications. In fact Seifert (\cite[Section 3.1.4]{Seifert97}) notes that 
``it turns out that in large regions of the interesting parameter space the shape of lowest energy is indeed axisymmetric for vesicles of spherical topology". We actually conjecture that \emph{for any given area and volume satisfying the isoperimetric inequality, and for any constant spontaneous curvature, the problem of minimizing \eqref{eq:helfrich} in the class of embedded, constrained surfaces has a solution, and it is axisymmetric}.
\medskip

After this work was completed, we became aware of the preprint \cite{Helmers}, which studies the $\Gamma$-limit of a diffuse-interface approximation of Helfrich's functional for two-phase axisymmetric surfaces, and where many of the same technical difficulties that we encounter are addressed. We are independently treating the sharp-interface case of two-phase axisymmetric surfaces in \cite{ChoMorVen}.

\subsection{The class of minimizers} An \textit{axisymmetric surface}, in the context of the present paper, is a surface $\S$ obtained by rotating a curve $\g$, contained in the $xz$ plane in $\R^3$, around the $z$-axis. Since, following the direct method of the calculus of variations, we want to find the minimum of $\H$ as the limit of a sequence of minimizers, we need the class of possible minimizers to be closed with respect to a reasonable topology. Simple curves alone are thus not sufficient, as a curve as in Figure \ref{fig1}-left  can be obtained as a uniform limit of smooth simple curves. The curves falling in the two classes introduced below are regular enough to allow for a definition of a generalized Helfrich energy, surface area, and enclosed volume, and at the same time are closed under the convergence induced by $\H$. 

\paragraph{Notation.} Let $\g:[a,b]\to \R^2$, $t\mapsto (\g_1(t),\g_2(t))$, be a plane curve of class $C^1$. Denote $\dot \g :=d \g/dt$. Let $(\g):=\g([a,b])=\{\g(t): t\in [a,b]\}$ be the trace of $\g$ and let $\ell(\g)$ be its length. We mostly parametrize $\g$ on the interval $[0,1]$ with constant speed $|\dot \g|=\ell(\g)$, in some cases, where specified, we use the arc-length parameterization $|\dot \g|\equiv 1$ on the interval $[0,\ell(\g)]$.  

\begin{defi}
\label{def:g0} 
A curve $\g :[0,1]\to \R^2$ belongs to the class \GO of \emph{curves generating a genus-0 surface with bounded weak curvature} if and only if 
\begin{align}
	&\g \in C^1((0,1);\R^2) \cap W^{2,2}_{\textup{loc}}((0,1);\R^2) \label{eq:condreg}\\
	&|\dot \g(t) |\equiv \ell(\g)\quad \forall\, t\in(0,1), \label{eq:condspeed}\\		
	&\g_1(0) =\g_1(1)=0,\quad \g_1(t) > 0\quad \forall\, t\in(0,1). \label{eq:condpos}
\end{align}	
\end{defi}

\begin{defi}
\label{def:g1} 
A curve $\g:[0,1]\to \R^2$ belongs to the class \GI of \emph{curves generating a genus-1 surface with bounded weak curvature}  if and only if 
\begin{align}
	& \g\in W^{2,2}\val,\label{eq:condreg1}\\
	&|\dot \g(t) |\equiv \ell(\g)\quad \forall\, t\in [0,1],\\			
	&\g(0) =\g(1),\quad	\dot \g(0) =\dot \g(1),\quad	\g_1(t) >0\quad \forall\,t\in [0,1]. \label{eq:condpos1}
\end{align}
\end{defi}

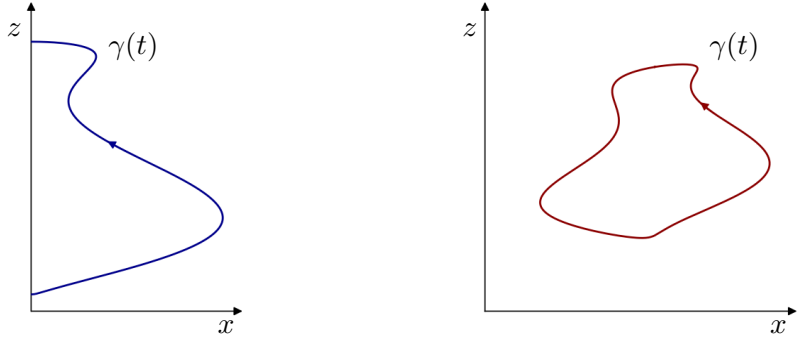
\begin{figure}[h]
\begin{minipage}{7cm}
	\labellist
		\hair 2pt
		\pinlabel $\g(t)$ at 130 320
		\pinlabel $z$ at -5 340	
		\pinlabel $x$ at 235 0	
	\endlabellist			
	\centering
	\includegraphics[height=4.5cm]{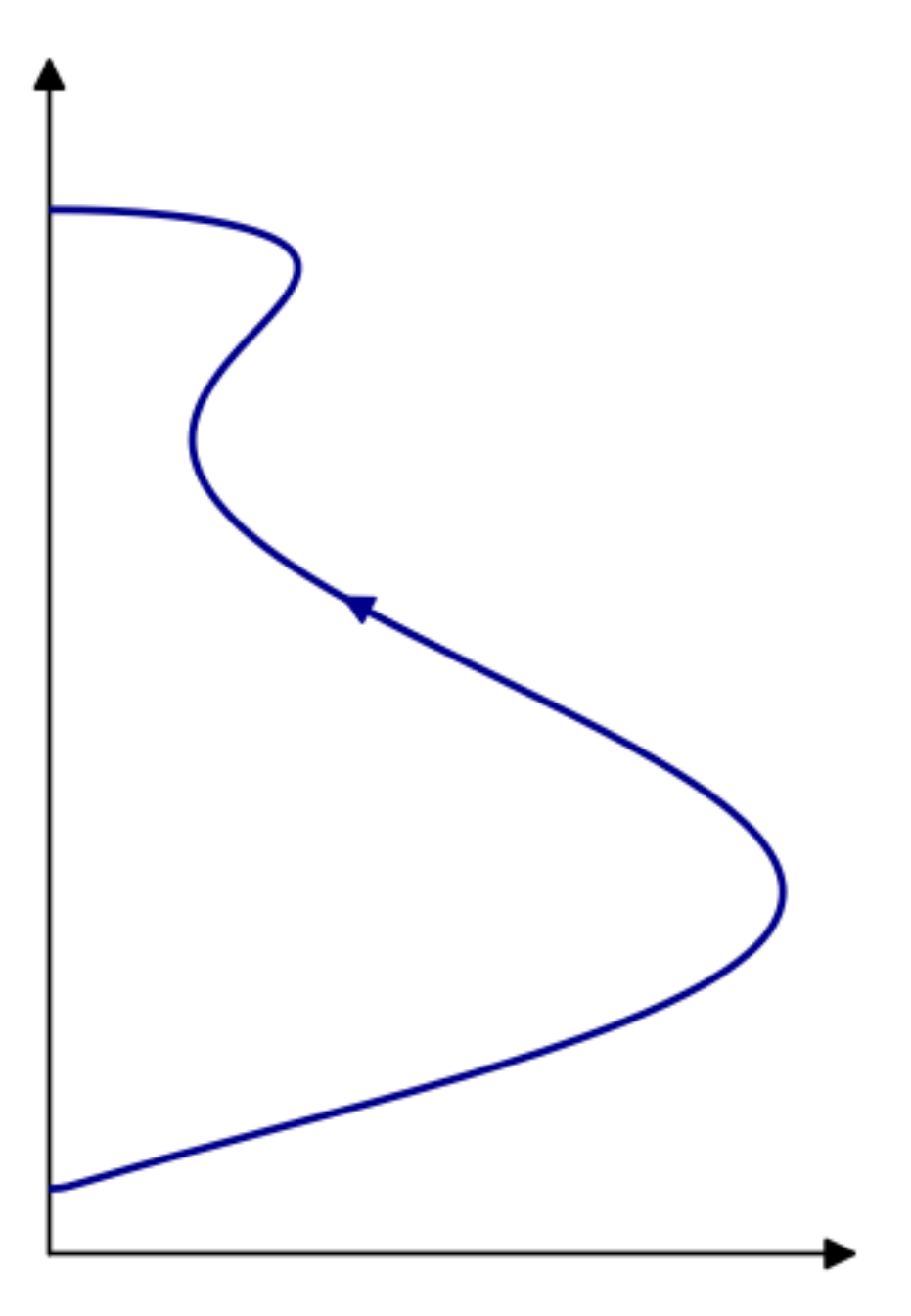}
\end{minipage}
\begin{minipage}{6cm}
	\centering
	\labellist
		\hair 2pt
		\pinlabel $\g(t)$ at 300 320
		\pinlabel $z$ at 0 340	
		\pinlabel $x$ at 350 0	
	\endlabellist				
	\includegraphics[height=4.5cm]{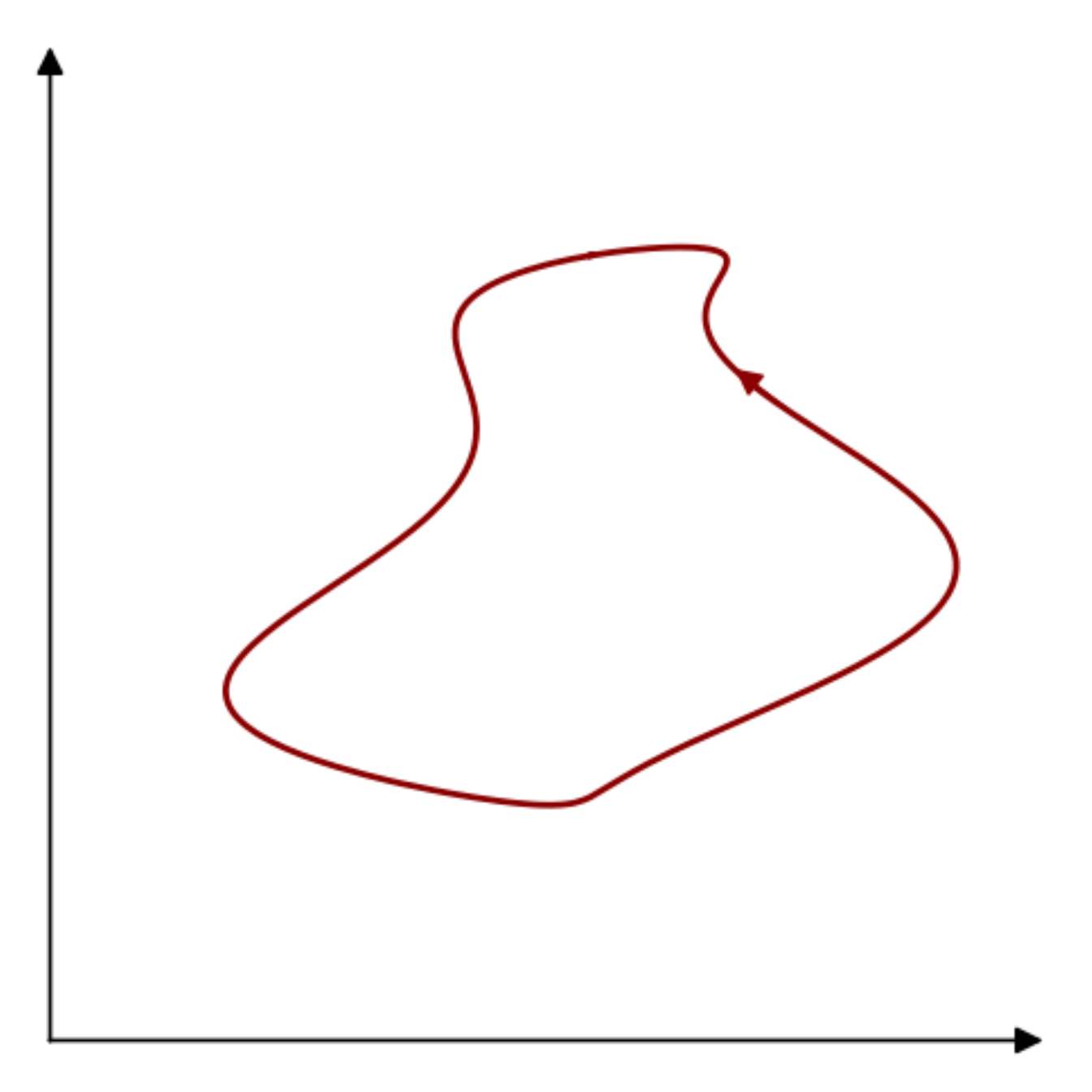}  
\end{minipage}
\caption{Generating curves in \GO (left) and \GI (right).}
\end{figure}
We note that the condition $\g_2(0)\neq\g_2(1)$ cannot be imposed in Definition \ref{def:g0}, since a curve with $\g_2(0) = \g_2(1)$ could be obtained as a continuous limit of curves in either \GO (see, e.g.,  Figure \ref{fig1}-left) or \GI.
\medskip
 
Let $\g$ be a curve as in \GO or \GI.  By rotating $\g$ around the $z$-axis we obtain the surface $\S$ parametrized by:
\begin{equation}
\label{def:parametr} 
	r(t,\theta)=\big[\g_1(t)\cos \theta,\ \g_1(t)\sin \theta,\ \g_2(t)\big],\qquad (t,\theta)\in [0,1]\times [0,2\pi].
\end{equation}	
If a surface $\S$ admits the parametrization \eqref{def:parametr}, we say that $\S$ \textit{is generated by} $\g$. 
\begin{figure}[h]
\begin{center}
\begin{minipage}{4.5cm}
\centering{
\labellist
		\hair 2pt
		\pinlabel $\S$ at 280 250
		\pinlabel $z$ at 185 315	
		\pinlabel $\g$ at 250 100
	\endlabellist				
\includegraphics[height=6cm]{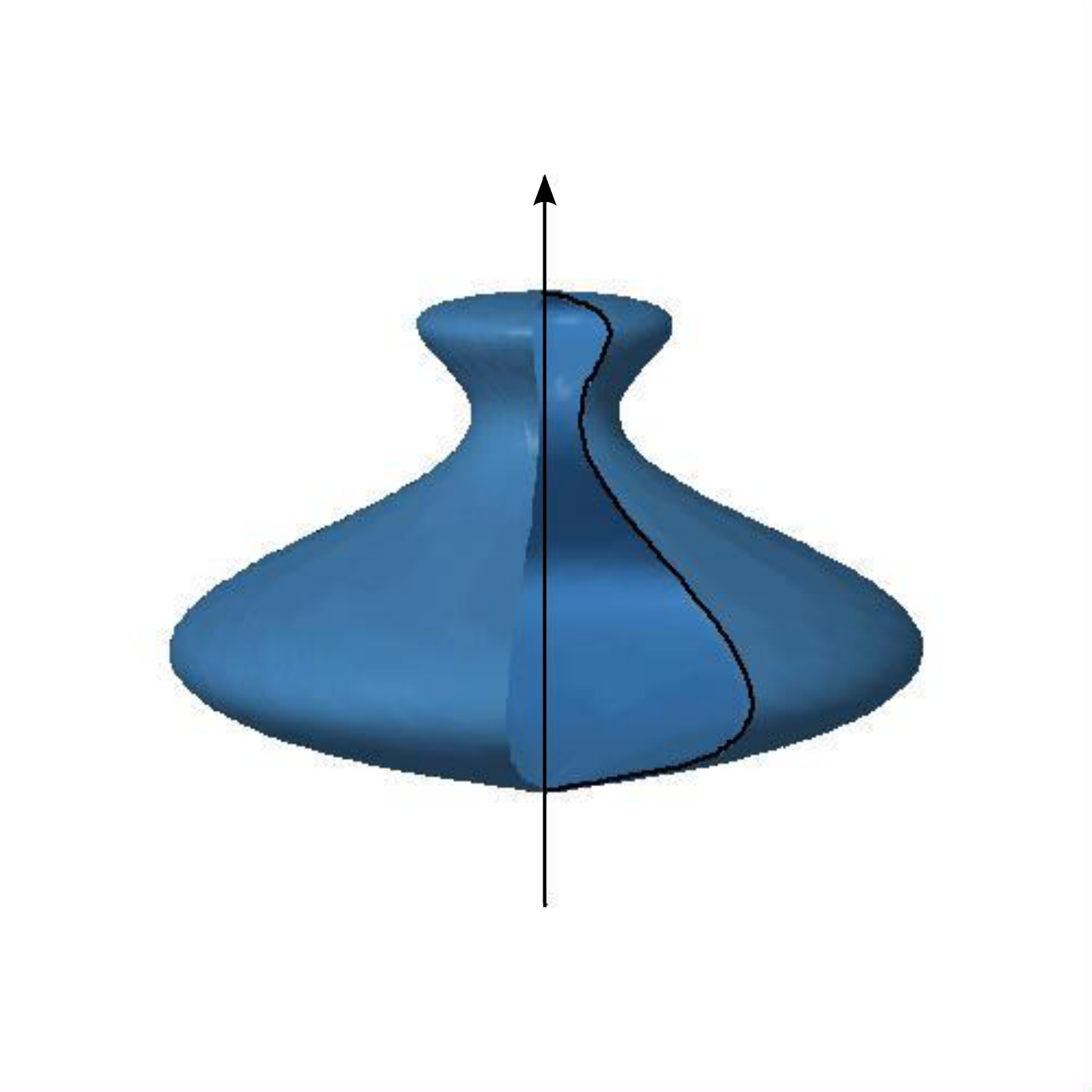}\\ 
\small    \vspace{-1cm} $\S$ generated by  $\g$ as in \GO}
\end{minipage}
\hspace{1.5cm}
\begin{minipage}{5cm}
\labellist
		\hair 2pt
		\pinlabel $\S$ at 385 250
		\pinlabel $z$ at 185 330	
		\pinlabel $\g$ at 345 120
	\endlabellist				
\centering{
\includegraphics[height=6cm]{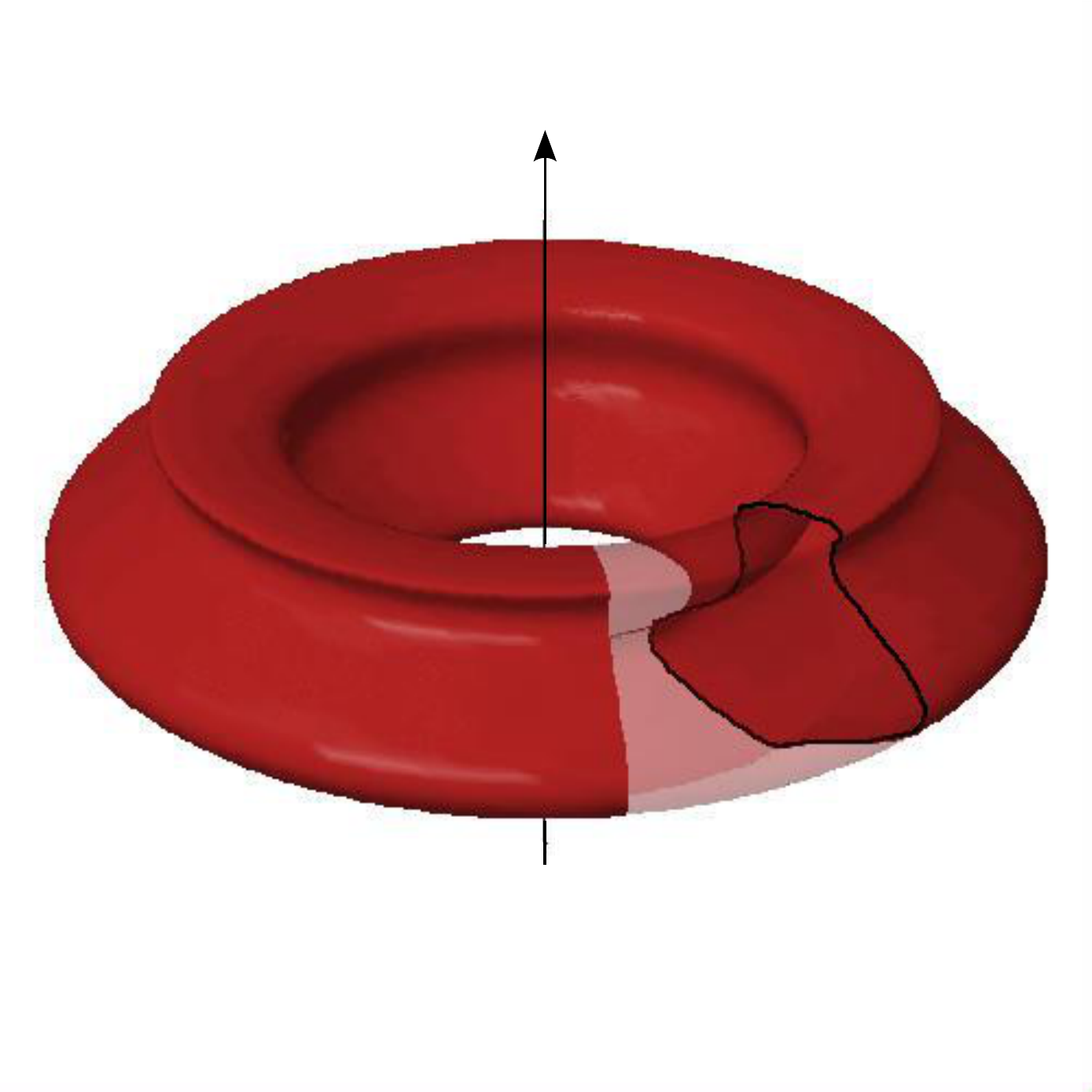}\\  
\small    \vspace{-1cm} \quad $\S$ generated by  $\g$ as in \GI}
\end{minipage}
\end{center}
\caption{Generated surfaces.}
\label{fig2}
\end{figure}
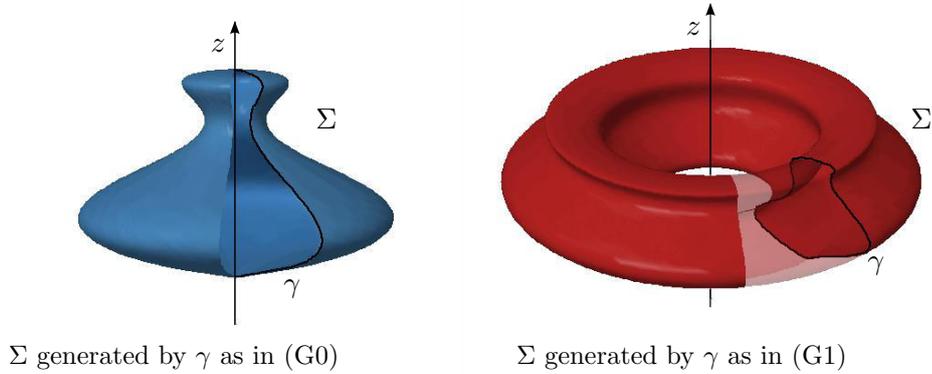

A standard computation (see Section \ref{ssec:surrev} below), shows that \emph{if} a curve $\g$ generates a smooth surface $\S$, the 2-dimensional surface area, the enclosed volume, and the principal curvatures of the generated surface are given by
\begin{align}
		 |\S| &=2\pi \int_0^1 \g_1|\dot \g|\, dt, & \vs &=\pi \int_0^1 \g_1^2\dot \g_2\, dt, \label{eq:firsdef} \\
		 k_1 &= \frac{(\ddot \g_2 \dot \g_1 -\ddot \g_1 \dot \g_2)}{|\dot \g|^3 },& k_2 &= \frac{\dot \g_2 }{\g_1|\dot \g| }. \label{eq:secdef}
\end{align}	
Since $H=k_1+ k_2$ and $K=k_1k_2$, the Helfrich energy can then be written as
\begin{align}
	\H(\S) &=   \int_\S \left\{  \frac{\kh}{2}(H-H_0)^2+\kg\, K\right\}dA \nonumber\\
		&= \int_0^1 \left\{ \frac{\kh}{2} \left(k_1+k_2  -H_0\right)^2 + \kg k_1 k_2\right\} 2\pi\g_1|\dot \g |\,dt. \label{eq:thirdef}
\end{align}
If the surface generated by $\g$ in \GO or \GI is not smooth, as is the case of Figure \ref{fig1}, we \emph{define} the \emph{generalized} 2-dimensional surface area, enclosed volume, principal curvatures and Helfrich energy of the generated surface by the quantities in \eqref{eq:firsdef}-\eqref{eq:secdef}. 

Since the integral of the Gaussian term of the energy is constant for a surface of fixed genus (see Section \ref{ssec:bonnet}), it is often disregarded in the analysis of minimizers. Nonetheless, since we are not imposing a fixed genus, nor a fixed number of components, we cannot drop this term. Furthermore, we note that the Gaussian term is expected to play an important role in the case of multiphase membranes \cite[page 1068]{BaumgartDWJ}.  
\medskip

The main result of this paper is the following.
\begin{theo}
\label{th:main}
Let $A,V>0$ be given such that
\begin{equation}
\label{eq:isop}
	V\leq \frac{A^{3/2}}{6 \sqrt{\pi}}.
\end{equation}
Assume that $\kh>0$, $\kg,H_0 \in \R$ such that $\frac{\kg}{\kh}\in (-2,0)$. Let $\A(A,V)$ denote the set of finite families $S=(\S_1,\ldots,\S_m)$, for some $m\in\N$ (not fixed), of axisymmetric surfaces generated by disjoint curves in \GO$\cup$ \textup{(G1)}, as in Definition \ref{def:g0} and Definition \ref{def:g1}, and satisfying the generalized area and volume constraints 
$$\sum_{i=1}^m|\S_i|=A,\qquad \sum_{i=1}^m \V(\S_i)=V.$$
Let $\H$ be the Helfrich energy functional defined in \eqref{eq:helfrich} and  let  
$$\ds \F: \mathcal A(A,V)\to \R \cup \{+\infty\},\qquad \F(S):= \sum_{i=1}^m \H(\S_i).$$
Then the problem 
$$ \min \left\{\F(S): S\in \A (A,V)\right\}$$
has a solution. 
\end{theo}

Condition \eqref{eq:isop} ensures that the constraints satisfy the isoperimetric inequality, so that the set $\mathcal A(A,V)$ is not empty. When \eqref{eq:isop} is an equality, the only element in $\mathcal A(A,V)$ is the sphere of area $A$, if it is a strict inequality, $\mathcal A(A,V)$ contains an infinite number of elements.

The range of the parameters $\kh$ and $\kg$ specified in the assumptions of Theorem \ref{th:main} is the mathematical range for which $\H$ is positive definite on the principal curvatures, or in other words, for which $\H(\S)$ controls the full squared norm of the second fundamental form of $\S$, which is absolutely crucial for any kind of analysis. On the other hand, the physical range in which these parameters are typically found is contained in the one we assume, see e.g. \cite{TemplerKhooSeddon} and \cite{BaumgartDWJ} (note that the latter cites the former, but inverting numerator and denominator, by mistake).

Note that the functional $\F$ does not depend on the reciprocal position of the components $\S_i$. Therefore, by translation along the vertical axis, we can transform a system with self intersections into one with the same energy and without crossings, thus avoiding unphysical situations.

\subsection{Discussion}
\label{ssec:discussion}

\begin{rem} 
\label{rem:index}
\textit{On the index $I(\g,p)$ of a curve.} Even if a curve has a smooth parametrization, it can generate a surface of revolution with singularities, which cannot represent any physical lipid bilayer (Figure \ref{fig1}). A way to restrict to physical surfaces is to prescribe the index of the system of generating curves.

If $\g$ is a closed curve, $p\in \R^2\backslash (\g)$, let $I(\g,p)$ be the index of $\g$ with respect to $p$ \cite[Chapter II, Section 1.8]{cartan}. If $\g$ is not closed, $\g_1(0)=\g_1(1)=0$ and $\g_1 \geq 0$, we can extend it symmetrically with respect to the $z$-axis in order to define its index. For a system of surfaces $S=(\S_1,\ldots,\S_m)$, generated by $(\g_1,\ldots,\g_m)$, and $p\in \R^2\backslash \cup_{i=1}^m(\g_i)$ define $I(S,p):=\sum_{i=1}^m I(\g_i,p).$ Note that  if $E\subset \R^2$ is a smooth connected bounded open set and $\g$ is counterclockwise parametrization of $\partial E$, then $E=\{p\in \R^2: I(\g,p)=1\}$ and $\R^2 \backslash \overline E = \{p\in \R^2: I(\g,p)=0\}$.   More generally,  points with index 1 represent the internal volume of a vesicle also for surfaces generated by curves which are not the parametrization of a boundary (as in Figure \ref{fig1}-left).  In order to eliminate situations like Figure \ref{fig1}-right, we can look for minimizers in the class of systems of surfaces $S \in \mathcal A(A,V)$ such that, additionally,  
\begin{equation}
\label{eq:index}
	I(S,p)\in \{0,1\}\ \text{for a.e. }p\in \R^2.
\end{equation}

The advantage in using the index as a condition, is its continuity with respect to uniform convergence of the curves, and thus its compatibility with the convergence induced by the bound $\F(S^n)\leq \Lambda$ (see Section \ref{ssec:plan} below and Definition \ref{def:systconv}), so that the proof of Theorem \ref{th:main} can be directly used to prove existence of minimizers in $\mathcal A(A,V)$ satisfying \eqref{eq:index}.  On the other hand, by removing the condition on the index, the minimizer of Helfrich functional could actually be found in non-embedded surfaces, like Delauney surfaces or the Wente torus.

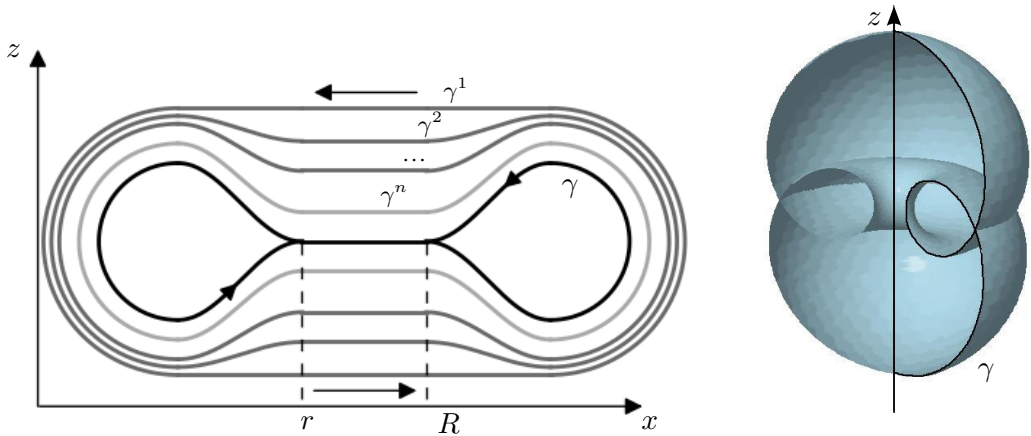
\begin{figure}[h!]
\begin{center}
\begin{minipage}{5cm}
	\labellist
	\hair 2pt
	\pinlabel {\footnotesize $\g^1$} at 247 190
	\pinlabel {\footnotesize $\g^2$} at 233 171
	\pinlabel {\footnotesize $\cdot\! \cdot\! \cdot$} at 224 151		
	\pinlabel {\footnotesize $\g^n$} at 213 131			
	\pinlabel $\g$ at 310 135
	\pinlabel $r$ at 165 5	
	\pinlabel $R$ at 243 5	
	\pinlabel $z$ at 3 212
	\pinlabel $x$ at 354 5		
	\endlabellist			
\centering
\noindent
\includegraphics[height=5.4cm]{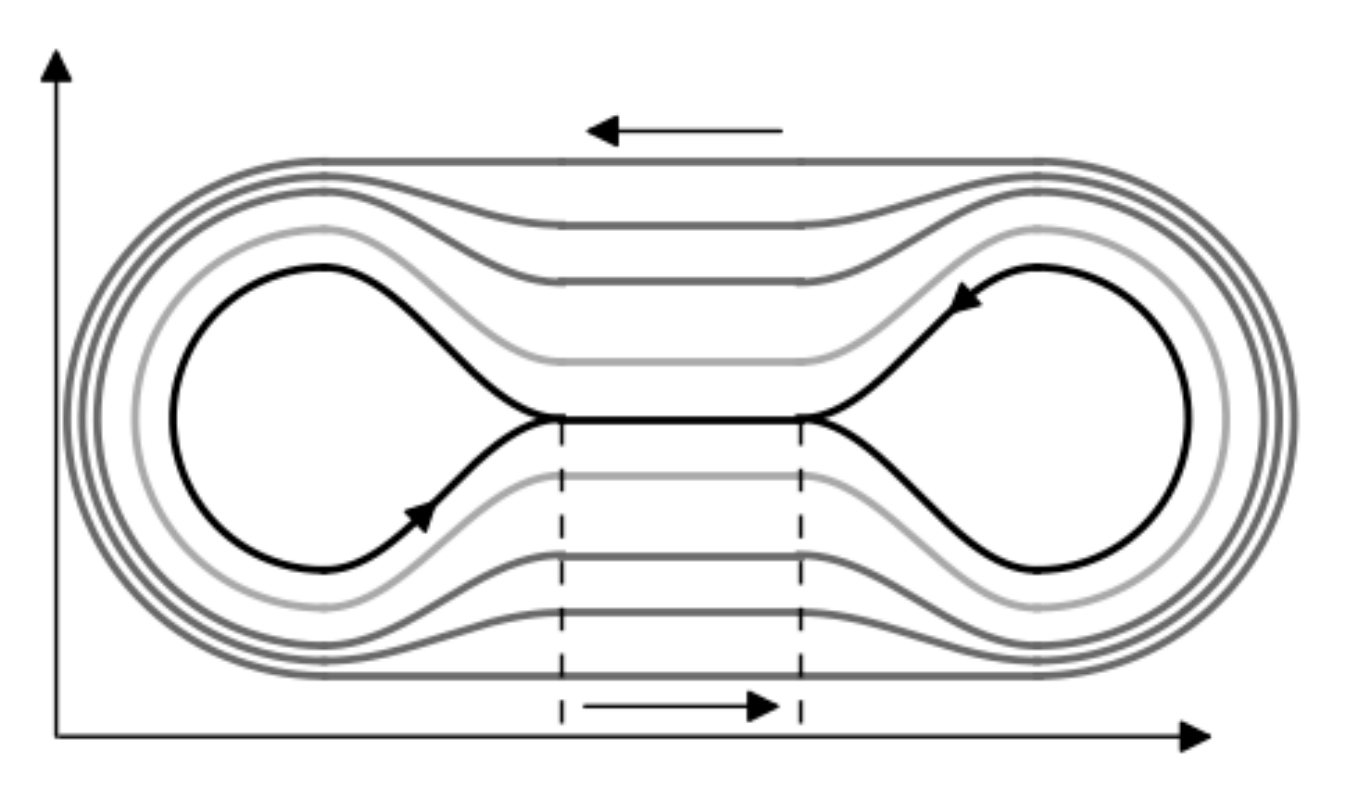}
\end{minipage}
\hspace{4cm}
\begin{minipage}{4cm}
\labellist
	\hair 2pt
	\pinlabel $\g$ at 270 40	
	\pinlabel $z$ at 160 390
\endlabellist			
\centering{
\includegraphics[height=5.4cm]{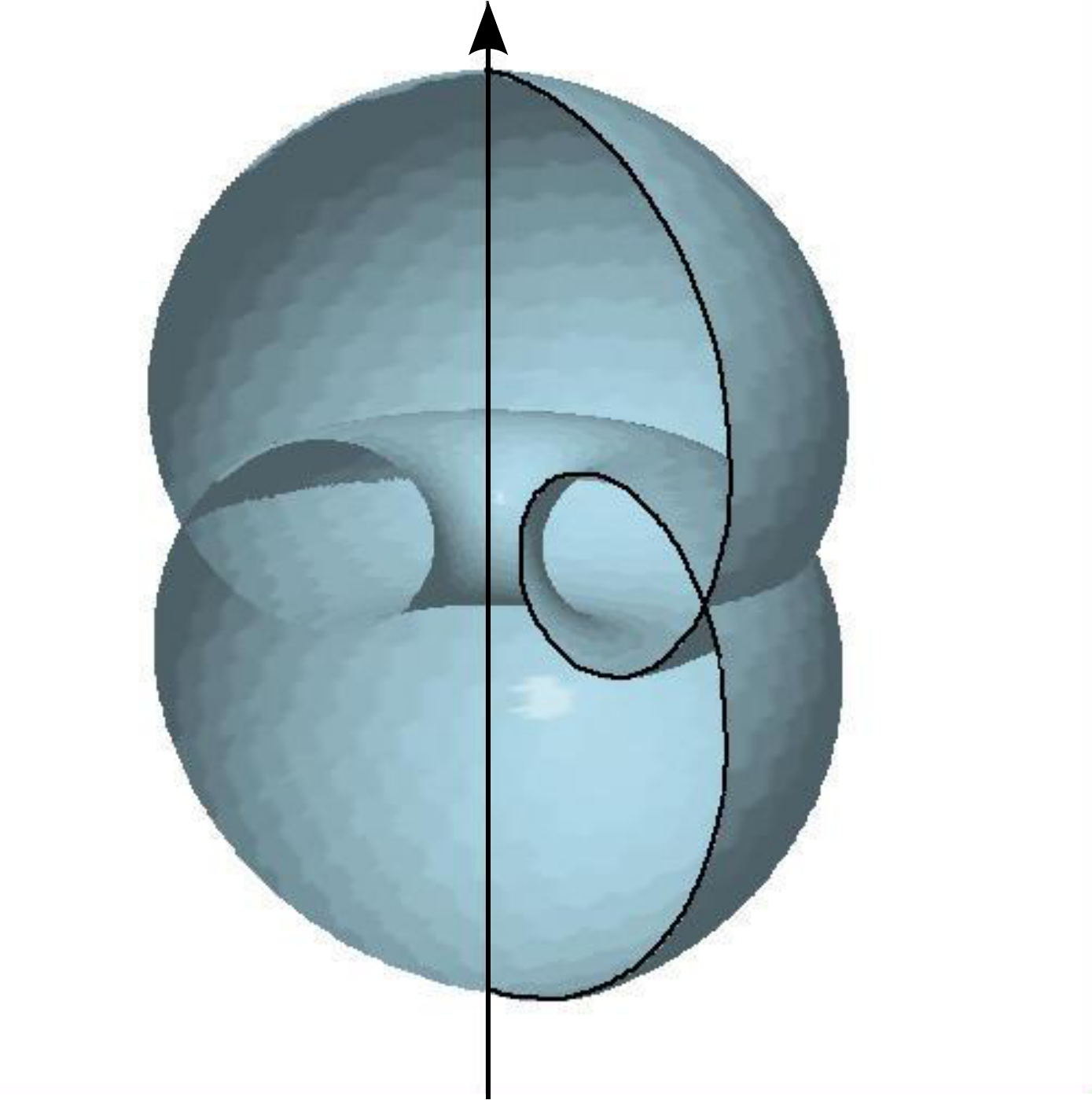} 
}
\end{minipage}
\end{center}
\vspace{-0.5cm}

\caption{smooth parametrizations of self-intersecting curves.}
\label{fig1}
\end{figure}

\end{rem}

\emph{On the lack of $W^{2,2}$-regularity.} In Corollary \ref{cor:reg} we prove that any curve that generates a surface with bounded Helfrich energy is $W^{2,2}$-regular on any stretch at positive distance from the $z$-axis. Since the area element vanishes on the $z$-axis, the second derivative of a curve in (G0) need not be square-integrable near the intersection with the $z$-axis, and therefore we cannot expect to control the  $L^2$-norm of $\ddot \g$. Loosely speaking, the reason why the regularity of a surface $\S$ does not imply the same regularity for the generating curve is simply the fact that a function (e.g., $|x|^{-1}$) can be integrable in a neighborhood of the origin in $\R^2$, but not in $\R$. For example, let $0<\delta<1$ and consider the curve defined by 
%
%
$$ \g_1(t):=\frac23 \left(1-(1-t)^{3/2}\right),\qquad \g_2(t):=\frac 23 t^{3/2},\qquad t\in [0,\delta].$$
\noindent Clearly, $\dot \g \in C^0([0,\delta];\R^2)$, $|\dot \g|\equiv 1$, $\ddot \g_2 \notin L^2((0,\delta);\R^2)$, and a quick computation shows that $ k_1^2\g_1 \sim 1/4$, $k_2^2 \g_1 \sim 1$, as $t\to 0^+$. Therefore, $k_1,k_2$ are square-integrable with respect to the area measure, on $(0,\delta)$. 
\medskip

\emph{On the generalized area and varifolds.} Returning to the example in Figure \ref{fig1}-left, we note that there could be two ways to describe the area of the middle annulus of the revolution surface generated by this curve. If we see it as a single layer of membrane, it should simply measure $4\pi(R^2-r^2)$. On the other hand, if we obtained this curve as a limit of a sequence of simple curves, where the vertical distance between the layers of the membrane collapsed to zero in the stretch between $r$ and $R$, it should be seen as a double layer, and measure $8\pi(R^2-r^2)$. Since we are imposing a constraint on the total area of the membrane, and the objects we wish to describe are \emph{closed} vesicles, the second interpretation, in which the horizontal stretch represents two overlapping layers, should be preferred. This is a drawback of modeling a three-dimensional object (i.e. a membrane with positive thickness) as a two-dimensional one: the only way we have to represent adjacent layers is to allow them to overlap, even though overlapping is not a physical possibility for the original three-dimensional membranes. 

Using parametrized curves to model membranes goes naturally in this second direction. The generalized area defined in \eqref{eq:firsdef}, in fact, following the parametrization, counts the multiplicity of every self-intersection. In this respect, we obtain the same result as if we described the curves with varifolds.  Without going into details (which can be found in \cite{Hutchinson86}), we may think of the weight-measure of an integral varifold $V$ associated to a 2-rectifiable set $\S\subset \R^3$ as $\mu:=\theta \mathcal H^2\llcorner \S$, where $\mathcal H^2$ is the two-dimensional Hausdorff measure, and $\theta$ is a measurable $\N$-valued function, which represents the \emph{density} of $V$. For example, if the surface $\S$ generated by $\g$ as in Figure \ref{fig1}-left is obtained as limit of smooth surfaces $\S^n$, we would have $\theta^n\equiv 1$, $\theta^n \mathcal H^2\llcorner \S^n \weakto \theta \mathcal H^2\llcorner \S$, in the sense of measures, $\theta(x) =1$ for $|x|\leq r$ and $|x|>R$, and $\theta(x)=2$ for $r<|x|<R$, counting the horizontal stretch twice, exactly as definition \eqref{eq:firsdef} does with curves.


\subsection{Plan of the paper and structure of the proof}
\label{ssec:plan}
In Section \ref{ssec:bonnet} we recall the main relations between mean curvature, Gaussian curvature and the Euler characteristic, and we show that Helfrich's energy controls the $L^2$-norm of the principal curvatures, with respect to the area measure. In Section \ref{ssec:surrev} we derive the main geometrical quantities for surfaces of revolution. In the following subsections we study which properties of a general generating curve can be obtained from the $L^2$ bound on the principal curvatures, we estimate the length of the curves (Section \ref{ssec:lengthbd}), we control the regularity of the tangents on the $z$-axis (Section \ref{sec:card}), and the total variation of $\dot \g_1$ (Section \ref{ssec:oscbd}). 

The proof of Theorem \ref{th:main} is given in Section \ref{sec:existence}, following the direct method of calculus of variations. The standard procedure consists in showing that sub-level sets of the functional $\F$ are compact, so that from a minimizing sequence $S^n$ it is possible to extract a converging subsequence $S^{n_k}\to S$, and then proving that $\F$ is lower-semicontinuous, so that $S$ is a global minimum for $\F$. 

The crucial ingredients, which are not prescribed by the general direct method and have to be chosen for each specific problem, are the set of minimizers, which we discuss at the beginning of Section \ref{sec:existence}, and the topology with respect to which compactness and lower-semicontinuity must be verified. Regarding the topology, it is natural to expect (or request) at least strong $W^{1,1}$ convergence for the generating curves, in order to preserve the surface area in the limit. Moreover, it is straightforward to show that the second fundamental form of $\S^n$ is uniformly bounded in $L^2$, but \textit{only with respect to the surface area measure}. Hence the necessity to study simultaneously convergence of the curvatures (expressed via $\g^n,\dot \g^n, \ddot \g^n$) and of the area measure $\mu_{\g^n}$. A suitable tool for this purpose is provided by the \emph{measure-function pairs} introduced in \cite{Hutchinson86}. The main definitions and theorems regarding measure-function pairs are recalled in Section \ref{ssec:mt}. The main body of the proof consists then in the lower-semicontinuity result (Proposition \ref{prop:lscS}) and in the compactness result (Proposition \ref{prop:compsurf}). Continuity of the constraints  follows from the choice of the topology, and it is described in Section \ref{ssec:constr}. All these steps are summarized in Section \ref{ssec:proof}, where the complete proof of Theorem \ref{th:main} is given. 

\section{Preliminaries and geometrical inequalities}
\label{sec:geoineq}
\setcounter{equation}{0}
\setcounter{theo}{0}

\subsection{Gauss-Bonnet theorem and positive definiteness of $\H$}
\label{ssec:bonnet}
Let $\chi(\S)$ be the Euler-Poincar\'e characteristic of a compact surface $\S$, see, e.g., \cite[Proposition 3, Section 4-5]{doCarmo}. Recall that every compact connected surface is homeomorphic to a sphere with a certain number $g$ of handles, and the number $ g=\frac{2-\chi(\S)}{2} $ is called the \textit{genus} of $\S$.
The Gauss-Bonnet Theorem states that if $\S$ has no boundary, then 
\begin{equation*}
	 \int_\S K\, dA = 2\pi \chi(\S).
\end{equation*}
Since we are dealing with surfaces of revolution, we are interested only in two cases:
\begin{itemize} 
	\item curves in (G0) (Figure \ref{fig2}-left), which generate surfaces homeomorphic to a sphere, hence $g=0$, $\chi(\S)=2$, and $\int K\, dA =4\pi$, and
	\item curves in (G1)  (Figure \ref{fig2}-right), which generate surfaces homeomorphic to a torus, hence $g=1$, $\chi(\S)=0$, and $\int K\,dA=0$.
\end{itemize}

The next Lemma contains the fundamental coercivity estimate for Helfrich's functional. It can be considered a standard observation (see e.g. \cite{BellettiniMugnaiTA}), but we report the proof for completeness.
\begin{lemma}
\label{lemma:totalbound}
Let $\S$ be generated by $\g\! \in\,$\GO\,$\cup$\,\GI, and let $\kh>0$, $\kg,H_0 \in \R$ such that $\frac{\kg}{\kh}\in (-2,0)$. Then there exists $C>0$ such that
\begin{equation}
\label{eq:totalbound}
	\int_\S (k_1^2 +k_2^2 )\,dA \leq C\big(|\S| +\H(\S)\big).
\end{equation}
\end{lemma}
\begin{proof}
Let $\lambda \in \R$, note that 
\begin{equation*}
	\frac12 (k_1 +k_2)^2 +\lambda k_1k_2 = \frac12 (k_1^2 +k_2^2) + (1+\lambda)k_1k_2 \geq \frac{1-|1+\lambda|}{2}(k_1^2+k_2^2),
\end{equation*}
and the coefficient in front of the last term is positive if and only if $\lambda \in (-2,0)$. For all $\e>0$ it holds
$$ 
	\frac{H^2}{2} = \frac{(H-H_0+H_0)^2}{2}\leq \frac {1+\e}2 (H-H_0)^2 + \frac{1+\e}{2\e}H_0^2,
$$
and thus
$$
	\frac {1+\e}2 (H-H_0)^2 + \frac{1+\e}{2\e}H_0^2 +\lambda(1+\e) K \geq \frac{1-|1+\lambda(1+\e)|}{2}(k_1^2+k_2^2).
$$
Choosing $\kh>0$, $\lambda = \kg/\kh \in (-2,0)$ and $\e>0$ such that $(1+\e)\kg/\kh \in (-2,0)$, we get
$$ 
	\frac{\kh}{2}(H-H_0)^2 + \kg K + c_1 H_0^2 \geq c_2(k_1^2+k_2^2),
$$
where $c_1=\kh/2\e$ and $c_2=\frac{\kh-|\kh+\kg(1+\e)|}{2(1+\e)}>0$. Integrating on $\S$ we obtain the thesis.
\end{proof}
We note that since we restrict to surfaces of revolution, the genus of which can only be 0 or 1, we could extend the range of parameters to $\kh/\kg >-2$. Indeed, if $\lambda =\kg/\kh \geq 0$, by Gauss-Bonnet theorem 
\begin{align*}
	\int_\S \frac12 (k_1 +k_2)^2 +\lambda k_1k_2\, dA&= \int_\S\frac12 (k_1^2 +k_2^2) + (1+\lambda)k_1k_2\, dA\\
		& \geq \frac12 \int_\S k_1^2+k_2^2\,dA + 2\pi \chi(\S).
\end{align*}
For a family of surfaces $S=(\S_1,\ldots,\S_n)$, we can then find a constant $C>0$ such that
$$ 	\sum_{i=1}^n\left(\int_{\S_i} (k_{1,i}^2 +k_{2,i}^2 )dA + \#\{\S_i \in S : g(\S_i)=0\} \right) \leq C\big(|\S| +\H(\S)\big). $$
 Since in physical applications the parameters are as in the assumptions of Lemma \ref{lemma:totalbound} (see, e.g., \cite{BaumgartDWJ}, \cite{TemplerKhooSeddon}), we will not use this estimate on the cardinality of the system, and rely only on \eqref{eq:totalbound}.

\subsection{Surfaces of revolution}
\label{ssec:surrev}

From this Section onwards, we restrict to surfaces of revolution. We start by deriving the geometrical quantities involved in the study of Helfrich's functional.

These computations can be found also, e.g., in \cite[Section 3-3, Example 4]{doCarmo} (with opposite orientation). With the parametrization 
\begin{equation*}
	r(t,\theta)=\big[\g_1(t)\cos \theta,\ \g_1(t)\sin \theta,\ \g_2(t)\big],\qquad (t,\theta)\in [0,1]\times [0,2\pi],
\end{equation*}	
we compute the tangent vectors
\begin{align*}
	r_t:=\frac{\partial}{\partial t}r(t,\theta)&= \big[\dot \g_1(t)\cos \theta,\ \dot\g_1(t)\sin \theta,\ \dot\g_2(t)\big], \\
	r_\theta:=\frac{\partial}{\partial \theta}r(t,\theta)&= \big[-\g_1(t)\sin \theta,\ \g_1(t)\cos \theta,\ 0\big].
\end{align*}
Note: $r_t\cdot r_\theta=0$, i.e., the tangents are always orthogonal. The first fundamental form is given by
$$	
g(t,\theta)=\left[
	\begin{array}{cc}
		E & F\\
		F & G 
	\end{array}	
	\right]=
	\left[
	\begin{array}{cc}
		r_t \cdot r_t & r_t \cdot r_\theta\\
		r_\theta \cdot r_t & r_\theta \cdot r_\theta 
	\end{array}	
	\right]
	=\left[
	\begin{array}{cc}
		|\dot \g(t)|^2 & 0 \\
 		0		& \g_1(t)^2
	\end{array}	
	\right],
$$	
$$ \sqrt{ g} := \sqrt{\det (g_{ij})}=\g_1(t)|\dot \g(t)|.$$
Note that the first fundamental form does not depend on the longitude parameter $\theta$. The normal vector can be oriented inwards or outwards, depending on the direction of $\g$. 
\begin{align*} 
	n(t,\theta) = \frac{r_t \times r_\theta}{\sqrt g}&=\frac{1}{\g_1(t)|\dot \g(t)|}\big[ -\g_1(t)\dot \g_2(t) \cos \theta,-\g_1(t)\dot \g_2(t) \sin \theta,\ \g_1(t)\dot \g_1(t)\big]\\
		&= \frac{1}{|\dot \g(t)|}\big[ -\dot \g_2(t) \cos \theta,-\dot \g_2(t) \sin \theta,\ \dot \g_1(t)\big].
\end{align*}		
For the computation of the second fundamental form  we make use of a constant-speed parametrization
\begin{align*}
	n_t:= \frac{\partial}{\partial t}n(t,\theta) &=\frac{1}{|\dot \g(t)|}  \big[ -\ddot \g_2(t) \cos \theta,-\ddot \g_2(t) \sin \theta,\ \ddot \g_1(t)\big]\\
	n_\theta:= \frac{\partial}{\partial \theta}n(t,\theta) &=\frac{1}{|\dot \g(t)|} \big[ \dot \g_2(t) \sin \theta,-\dot \g_2(t) \cos \theta,\ 0\big]
\end{align*}
$$	II(t,\theta)=\left[
	\begin{array}{cc}
		L & M\\
		M & N 
	\end{array}	
	\right]=-\left[
	\begin{array}{cc}
		n_t \cdot r_t & n_t \cdot r_\theta\\
		n_\theta \cdot r_t & n_\theta \cdot r_\theta 
	\end{array}	
	\right]
	=\frac{1}{|\dot \g|} \left[
	\begin{array}{cc}
		\ddot \g_2 \dot \g_1 -\ddot \g_1 \dot \g_2  & 0 \\
 		0		& \g_1 \dot \g_2
	\end{array}	
	\right].
$$	
We  can then express the Gaussian curvature 
$$ K=k_1 k_2 = \frac{LN -M^2}{EG - F^2}=\frac{(\g_1 \dot \g_2)(\ddot \g_2 \dot \g_1 -\ddot \g_1 \dot \g_2)}{(\g_1)^2|\dot \g|^4 },$$
 (twice) the mean curvature
\begin{equation*} 
	H=k_1 +k_2 = \frac{LG-2MF +NE }{EG - F^2}=\frac{\g_1^2 (\ddot \g_2 \dot \g_1 -\ddot \g_1 \dot \g_2) + \g_1\dot \g_2 |\dot \g|^2}{(\g_1)^2|\dot \g|^3 },
\end{equation*}
and the principal curvatures
\begin{align*}
	k_1&=\frac{\ddot \g_2 \dot \g_1 -\ddot \g_1 \dot \g_2}{|\dot \g|^3}		\quad \text{(meridian),}&&
	k_2=\frac{\dot \g_2}{\g_1|\dot \g| }\quad \text{(parallel).}
\end{align*}	
Note that $k_1$ is just the curvature of $\g$, with the sign depending on the orientation. Let $\L^1\llcorner_{[0,1]}$ be the one-dimensional Lebesgue measure restricted to the closed interval $[0,1]$, define the Radon Measure 
\begin{equation}
\label{eq:defmug}
	\mu_\g:= 2\pi \g_1|\dot \g|\L^1\llcorner_{[0,1]}.
\end{equation}
 The area of the generated surface is given by
\begin{equation}
\label{eq:defarea}
	|\S|=\int_\S dA=\int_0^{2\pi}\int_0^1 \sqrt{g(t)}\,dt\,ds= 2\pi \int_0^1 \g_1(t)\, |\dot \g(t)|\, dt= \int_0^1 d\mu_\g,
\end{equation}
and the enclosed volume by
$$ 
	\vs=\pi \int_0^1 \g_1(t)^2\dot \g_2(t)\,  dt.
$$
Owing to \eqref{eq:totalbound}, and recalling that we are using a constant-speed parametrization,
\begin{align*}
	C(\H(\S) +|\S|) \geq \int_\S k_1^2 + k_2^2 \, dA = 2\pi \int_0^1 \left( \frac{|\ddot \g|^2}{|\dot \g|^4} + \frac{\dot \g_2^2}{\g_1^2|\dot \g|^2} \right)\g_1|\dot \g|\, dt.
\end{align*}
Since $\H$ and the total area are invariant under reparametrizations of $\g$, we can write the last inequality in the case of arc-length parametrization, obtaining the crucial estimate 
\begin{equation*}
	C(\H(\S) +|\S|) \geq   \int_0^{\ell(\g)} \left( |\ddot \g|^2\g_1 + \frac{\dot \g_2^2}{\g_1} \, \right) dt.
\end{equation*}

\subsection{A bound on the length}
\label{ssec:lengthbd}
In this subsection we derive a uniform bound on the length of the generating curves. 
\begin{lemma} 
\label{lemma:lengthbound}
Let $\g$ be curve generating a revolution surface $\S$ as in \GO or \GI
. Then 
$$
	\frac{|\S|}{2\pi\, \dms}\leq \ell(\g) \leq \frac{\sqrt{ |\S|}}{2\pi}\left\{ \left( \int_\S k_1^2\, dA \right)^{1/2} + \left( \int_\S k_2^2\, dA \right)^{1/2}  \right\}.
$$	
\end{lemma}
\begin{proof}
In order to obtain the left inequality, 
we compute
$$
	|\S| = 2\pi\int_0^1 \g_1(t)|\dot \g(t)|\, dt\leq 2\pi\,\ell(\g)\max_{t \in [0,1]}|\g_1(t)|\leq 2\pi\,\ell(\g)\dms.
$$
Regarding the right inequality, it is not restrictive to assume that $\g$ is parametrized by arc-length, on the interval $[0,L]$, where $L:=\ell(\g)$, so that 
\begin{align}
	|\dot \g| &\equiv 1,\label{eq:dotgarcl}\\
	| k_1 | &= \frac{|\ddot \g \times \dot \g|}{|\dot \g|^3}= |\ddot \g|. \label{eq:k1arcl}
\end{align}
We compute
\begin{equation}
\label{eq:Lbd1}
	L = \int_0^L 1\, dt = \int_0^L |\dot \g(t)|^2dt = \int_0^L \dot \g_1(t)\dot \g_1(t) + \dot \g_2(t) \dot \g_2(t)\, dt.
\end{equation}
By either (G0) or (G1) above, $\dot\g_1\g_1(L)=\dot\g_1\g_1(0)$, and $\g_1(t)> 0$ for all $t\in(0,L)$, so the first term of the right-hand side of \eqref{eq:Lbd1} becomes
\begin{align*}
	\int_0^L \dot \g_1(t)\dot \g_1(t)\, dt &= -\int_0^L  \g_1(t)\ddot \g_1(t)\,dt +\Big[\g_1(t)\dot\g_1(t)\Big]_{t=0}^{t=L} \\
		&=\int_0^L  \g_1(t)\ddot \g_1(t)\,dt \\
		&=\int_0^L  \sqrt{\g_1(t)}\left(\sqrt{\g_1(t)}\ddot \g_1(t)\right)\,dt \\
		&\leq \left( \int_0^L \g_1(t)\,dt \right)^{1/2}\left( \int_0^L \g_1(t)\ddot \g_1^2(t)\,dt \right)^{1/2}\\
		&\stackrel{\eqref{eq:defarea}}{=} \left( \frac{|\S|}{2\pi} \right)^{1/2}\left(\frac{1}{2\pi} \int_0^L  |\ddot \g(t)|^2\, 2\pi\g_1(t)\,dt \right)^{1/2}\\
		&\stackrel{\eqref{eq:k1arcl}}{=} \frac{\sqrt{ |\S|}}{2\pi}\left( \int_\S k_1^2\, dA \right)^{1/2}.
\end{align*}
On the other hand, the last term of \eqref{eq:Lbd1} gives
\begin{align*}
	\int_0^L \dot \g_2(t)\dot \g_2(t)\, dt & \leq \max_{t\in [0,L]}|\dot \g_2(t)| \int_0^L |\dot \g_2(t)|\, dt\\
		&\leq |\dot \g| \int_0^L \sqrt{\g_1(t)}\frac{|\dot \g_2(t)|}{\sqrt{\g_1(t)}}\, dt\\
		&\stackrel{\eqref{eq:dotgarcl}}{\leq} \left( \int_0^L \g_1(t)\,dt \right)^{1/2}\left( \int_0^L \frac{\dot \g_2^2(t)}{\g_1(t)}\,dt \right)^{1/2}\\
		&= \frac{1}{2\pi}\left( |\S| \int_\S k_2^2\, dA \right)^{1/2}.	
\end{align*}
\end{proof}
Noting that 
\begin{align*} 
	2\sqrt{ |\S|}\left\{ \left( \int_\S k_1^2\, \right)^{1/2} \hspace{-0.2cm}+ \left( \int_\S k_2^2\, \right)^{1/2}  \right\} &\leq |\S| + \left\{ \left( \int_\S k_1^2\, \right)^{1/2} \hspace{-0.2cm} + \left( \int_\S k_2^2\, \right)^{1/2}  \right\}^2 \\
	&\leq |\S| +  2\int_\S k_1^2 + k_2^2\, 
\end{align*}	
we get the following bound for a system of curves.
%
\begin{cor}
\label{cor:realbd}
Let $\g_i$, $i=1,\ldots,m$ be a finite family of curves in \GO or \GI generating the revolution surface $\S_i$.
Then
\begin{equation}
\label{eq:totalbounds}
	2\pi \sum_{i=1}^m\ell(\g_i) \leq \sum_{i =1}^m\left(\frac{|\S_i|}{2} +\int_{\S_i} k_{1,i}^2+k_{2,i}^2\,dA_i\right)	.
\end{equation}
\end{cor}

\subsection{Regularity of generators}
\label{ssec:reg}
\begin{defi}
\label{def:gengen}
We say that $\g:[0,1]\to\R^3$ is a \emph{generalized generator} if $\g$ is Lipschitz-continuous, $|\dot \g(t)|\equiv \ell(\g)$ and $\g_1(t) > 0$ for almost every $t \in (0,1)$, $\ddot \g \in L^1_{\textup{loc}}(\{\g_1>0\};\R^2)$ and 
\begin{equation}
\label{eq:L1bound}
	\int_0^1 k_1^2 +k_2^2\, d\mu_\g< C. 
\end{equation}
In particular, \eqref{eq:L1bound} implies that $\g \in L^2(\mu_{\g};\R^2)$.
\end{defi}

\label{sec:card}
\begin{lemma}[Internal regularity]
\label{lemma:w2loc}
Let $\g$ be as in Definition \ref{def:gengen}. For every subinterval $[a,b]\subset [0,1]\cap \{\g_1>0\}$ 
$$ \g\in W^{2,2}((a,b);\R^2),\quad \text{ and }\dot \g \text{ has a unique extension to } C^0([a,b];\R^2).$$
\end{lemma}
\begin{proof}
It holds	
\begin{equation}
\label{eq:strongw2}
	 \int_a^b k_1^2\,d\mu_\g = \int_a^b \frac{|\ddot \g|^2}{|\dot \g|^4} 2\pi |\dot \g|\g_1\, dt \geq \frac{2\pi}{ \ell(\g)^3}\min_{s\in[a,b]} \{\g_1(s)\} \int_a^b|\ddot \g|^2\, dt.
\end{equation}
Thus, $\ddot \g \in L^2((a,b);\R^2)$ and $\g\in W^{2,2}((a,b);\R^2)$. By standard Sobolev inclusions, $\dot \g \in W^{1,2}((a,b);\R^2)\hookrightarrow C^0((a,b);\R^2)$, and there is a unique function which extends $\dot \g$ to $C^0([a,b];\R^2)$. We denote this extension by $\dot \g$.
\end{proof} 
In particular, if $\g\in$(G1), then $\g\in W^{2,2}((0,1);\R^2)$ and $\dot \g$ has a unique extension to $ C^0([0,1];\R^2)$.


\begin{lemma}[Tangents on the $z$-axis]  
\label{lemma:tg}
Let $\g$ be as in Definition \ref{def:gengen}. Let $a,b\in [0,1]$ be such that $\g_1(a)=\g_1(b)=0$, $\g_1(t)>0$ for all $t\in  (a,b)$. Then, the limits of $\dot \g_2(t)$ as $t \to a^+$ and $t \to b^-$ exist, and
\begin{equation*}
	\lim_{t \to a^+}\dot \g_2(t) = \lim_{t \to b^-}\dot \g_2(t) = 0.
\end{equation*}
Moreover, either
$$ 
	\lim_{t\to a^+} \dot \g_1(t)=\ell(\g),\qquad 	\lim_{t\to b^-} \dot \g_1(t)=-\ell(\g),
$$
or
$$ 
	\lim_{t\to a^+} \dot \g_1(t)=-\ell(\g),\qquad 	\lim_{t\to b^-} \dot \g_1(t)=\ell(\g).
$$
\end{lemma}
\begin{proof}
It is not restrictive to assume $|\dot \g| \equiv 1$ as before. For any $(s,r)\in (a,b)$
\begin{align}
	\left |\dot \g_2^2(r)-\dot\g_2^2(s) \right| &=\left| \int_s^r 2 \dot \g_2 \, \ddot \g_2\, d\s\right| = 2\left| \int_s^r  \frac{\dot \g_2}{\sqrt{\g_1}} \, \ddot \g_2\sqrt{\g_1 }\, d\s\right| \nonumber\\
	 &\leq \int_s^r  \frac{\dot \g_2^2}{\g_1 } + |\ddot \g|^2 \g_1\, d\s \leq \frac{1}{2\pi}\int_s^r k_2^2 +k_1^2\, d\mu_\g \label{eq:cauchy}
\end{align}	

Since $\dot \g_2^2/\g_1$ and $|\ddot \g|^2 \g_1$ belong to $L^1(a,b)$ by \eqref{eq:L1bound}, we can define an absolutely continuous function
$$ G(r):= \int_a^r  \frac{\dot \g_2^2}{\g_1 } + |\ddot \g|^2 \g_1\, d\s, $$
satisfying 
\begin{equation}
\label{eq:Gt}
	\lim_{r\to a^+}G(r)=0.
\end{equation}
By \eqref{eq:cauchy} and \eqref{eq:Gt}, for every $\e>0$ there exists a $\delta>0$ such that
$$ \sup_{r,s \in (a,a+\delta)}\big| \dot \g_2^2(r)- \dot\g_2^2(s) \big|\leq \e,$$
thus, the limit of $\dot \g_2^2(s)$ as $s\to a^+$ exists. We can now prove that this limit is 0. Recall that $|\dot \g|=1$ and $\g_1(a)=0$ by hypothesis, then for all $\e>0$ we have
$$ \max_{s\in [a,a+\e]} \g_1(s)= \max_{s\in [a,a+\e]} \int_a^s \dot \g_1(\s)\, d\s\leq \e|\dot \g|=\e,$$
and
$$ 0\stackrel{\eqref{eq:Gt}}{=}\lim_{\e \to 0^+}\int_t^{t+\e} \frac{\dot \g_2^2}{\g_1}\, d\s \geq \limsup_{\e \to 0^+} \frac 1{\e} \int_t^{t+\e} \dot \g_2^2\, d\s =  \limsup_{\e \to 0^+} \dashint_t^{t+\e} \dot \g_2^2\, d\s.$$
Since the integrand is nonnegative, we conclude that $\lim_{t \to a^+}\dot \g_2(t) =0.$ The proof of the corresponding statement for the limit as $t \to b^-$ is identical. The statement on the limit of $\dot \g_1$ follows by the assumption $|\dot \g|^2= \ell(\g)^2=\dot \g_1^2+\dot \g_2^2$ and by the continuity of $\dot \g$ in $(a,b)$ obtained in Lemma \ref{lemma:w2loc}.
\end{proof}

\begin{cor}[Regularity]
\label{cor:reg}
Under the assumptions of Lemma \ref{lemma:tg}
$$
	\g \in W^{2,2}_{\textup{loc}}\big((a,b);\R^2\big),\quad \dot \g \in C^0([a,b];\R^2).
$$
\end{cor}
\begin{proof}
Since $\g_1>0$ in $(a,b)$, by Lemma \ref{lemma:w2loc} it holds $\g \in W^{2,2}_{\textup{loc}}\big((a,b);\R^2\big)$. By Lemma \ref{lemma:tg},  $\dot \g$ has a continuous extension to $\in C^0([0,1])$. 

\end{proof}

\subsection{A bound on the oscillations}
\label{ssec:oscbd}
\begin{lemma}		
\label{lemma:oscbd} 
Let $\g$ be as in Definition \ref{def:gengen}. Let $(a,b)\subseteq (0,1)$ be such that $\g_1(t)>0$ for all $t\in (a,b)$, then
\begin{equation}
\label{eq:oscillbd}
 	|\dot \g|\int_a^b k_1^2+ k_2^2\, d\mu_\g \geq 4\pi|\dot \g_1(b)-\dot \g_1(a)|
\end{equation}
and
\begin{equation}
\label{eq:oscillbd2}
 	\sqrt{|\dot \g|}\left(\int_a^b \k_1^2\, d\mu_\g + \int_a^b \frac{\dot \g_1^2(t)}{\g_1(t)}\, dt\right) \geq 2\sqrt{2\pi}|\dot \g_2(b)-\dot \g_2(a)|
\end{equation}
\end{lemma}

\begin{proof}
By Lemma \ref{lemma:w2loc}, or Corollary \ref{cor:reg}, we can assume that $\dot \g\in C^0([a,b];\R^2)$ and  $\g\in W^{2,2}_{\textup{loc}}((a,b);\R^2)$. Recall that  $|\dot \g|^2 =\dot \g_1^2 +\dot \g_2^2 \equiv \ell^2(\g)$ and $|k_1|=|\ddot \g|/\ell^2(\g)$. On $\{\dot \g_2\neq 0\}$ we compute
$$
	\dot \g_2^2=\ell^2(\g)-\dot \g_1^2,\qquad \dot \g_2 =\pm \sqrt{\ell^2(\g)-\dot \g_1^2}, \qquad \ddot \g_2 = \mp\frac{\dot \g_1\ddot \g_1}{\sqrt{\ell^2(\g)-\dot \g_1^2}},\qquad \ddot \g_2^2 =\frac{\dot \g_1^2\ddot \g_1^2}{\ell^2(\g)-\dot \g_1^2}\, .
$$	
Defining
$$
	\eta(t):=\left\{ 
		\begin{array}{l}
			1 \quad \text{if }\dot \g_2(t)\neq 0\\
			0 \quad \text{if }\dot \g_2(t) = 0,
		\end{array} 
	\right.
$$
we have
\begin{align}
	 \frac{\ell(\g)}{2\pi}\int_a^b \k_1^2\, dA &= \frac{\ell(\g)}{2\pi}\int_a^b \frac{|\ddot \g|^2}{\ell^4(\g)} 2\pi \g_1\ell(\g)\, dt =\frac{1}{\ell^2(\g)} \int_a^b (\ddot \g_1^2 +\ddot \g_2^2)\g_1\, dt\nonumber\\
	 & \geq \frac{1}{\ell^2(\g)} \int_a^b \eta \left( \ddot \g_1^2 + \frac{\dot\g_1^2\ddot \g_1^2}{\ell^2(\g)-\dot \g_1^2}\right)\g_1\, dt 
	 =  \frac{1}{\ell^2(\g)} \int_a^b  \eta \ddot \g_1^2\left( 1 + \frac{\dot\g_1^2}{\dot \g_2^2}\right)\g_1\, dt \nonumber\\
		&   = \frac{1}{\ell^2(\g)}  \int_a^b  \eta \ddot \g_1^2\left( \dot \g_2^2+\dot \g_1^2\right)\frac{\g_1}{\dot \g_2^2}\, dt = \int_a^b  \eta \ddot \g_1^2\, \frac{\g_1}{\dot \g_2^2}\, dt \nonumber\\
		&  \geq \inf \left\{    \int_a^b  \dot \phi^2\, \frac{\g_1}{\dot \g_2^2}\, dt : \phi \in W^{1,2}_{\textup{loc}}((a,b);\R^2),\ \phi(a)=\dot \g_1(a),\ \right.\label{eq:minimize}\\
		&\hspace{4cm} \left. \phi(b)=\dot \g_1(b),\ \dot \phi^2\, \frac{\g_1}{\dot \g_2^2}=0 \text{ on } \{\dot \g_2=0\}\right\}.\nonumber 
\end{align}
Let $\psi:= \g_1 / \dot \g_2^2$, $x_a:=\dot \g_1(a),\ x_b:=\dot \g_1(b)$. The unique minimizer of the above problem is the solution of the Euler-Lagrange equation
$$ \frac{d}{dt}(\dot \phi \psi)=0,\qquad \phi(a)=x_a,\quad \phi(b)=x_b,\quad \dot \phi^2\, \psi=0 \text{ on } \{\dot \g_2=0\}.$$
By integration, we compute
$$ \dot \phi = \frac{C}{\psi},\qquad \phi(t)=C'+C\int_a^t \frac{1}{\psi(t)}dt.$$
Define $J:= \int_a^b 1/\psi(t)\, dt$, imposing the boundary conditions we get
\begin{equation}
\label{eq:defphi}
	C'=x_a,\qquad C=\frac{x_b-x_a}{J},\qquad \phi(t)=x_a+\frac{x_b-x_a}{J}\int_a^t \frac{1}{\psi(t)}dt.  
\end{equation}
Note that if $\dot \g_2 \equiv 0$ in $(a,b)$, then $\dot \g_1$ is constant in $(a,b)$, and \eqref{eq:oscillbd} is trivially satisfied. If $\dot \g_2^2(t)>0$ in a point $t$, by continuity it is positive in an open interval containing t, and therefore $J>0$. Then, $\phi$ in \eqref{eq:defphi} is well-defined and, in particular $\dot \phi^2\psi=0$ on $\{\dot \g_2=0\}$. The minimum value is then given by
$$ \int_a^b \dot \phi^2(t)\, \psi(t)\, dt = \int_a^b \left( \frac{x_b-x_a}{J\, \psi(t)}\right)^2\psi(t)\, dt= \left( \frac{x_b-x_a}{J}\right)^2\int_a^b \frac{1}{ \psi(t)}\, dt = \frac{(x_b-x_a)^2}{J}.$$
Since
$$ J=\int_a^b \frac{1}{\psi(t)}\, dt = \ell(\g)\int_a^b \frac{\dot \g_2^2(t)}{\g_1(t)\ell(\g)}\, dt= \frac{\ell(\g)}{2\pi}\int_a^b \k_2^2\, dA,$$
by inserting the minimum value in \eqref{eq:minimize} and multiplying by $J$ we obtain
\begin{equation*}
 	\left(\int_a^b \k_1^2\, dA\right)\left(\int_a^b \k_2^2\, dA\right) \geq \frac{4\pi^2}{\ell^2(\g)}|\dot \g_1(b)-\dot \g_1(a)|^2.
\end{equation*}
Noting that
$$
	\left(\int_a^b \k_1^2\, dA\right)\left(\int_a^b \k_2^2\, dA\right)\leq \frac 14 \left( \int_a^b \k_1^2\, dA + \int_a^b \k_2^2\, dA \right)^2= \left(\frac 12 \int_a^b \k_1^2+ \k_2^2\, dA \right)^2
$$	
and taking the square root, we obtain \eqref{eq:oscillbd}.

In order to prove \eqref{eq:oscillbd2}, we follow the same computations, inverting the role of $\dot \g_1$ and $\dot \g_2$. Let
$$
	\eta(t):=\left\{ 
		\begin{array}{l}
			1 \quad \text{if }\dot \g_1(t)\neq 0\\
			0 \quad \text{if }\dot \g_1(t) = 0,
		\end{array} 
	\right.
$$
we have
\begin{align}
	 \frac{\ell}{2\pi}\int_a^b \k_1^2\, d\mu_\g &= \frac{1}{\ell^2}\int_a^b (\ddot \g_1^2 +\ddot \g_2^2)\g_1\, dt\nonumber\\
	 &\geq \frac{1}{\ell^2}\int_a^b \eta \left(  \frac{\dot\g_2^2\ddot \g_2^2 }{\ell^2-\dot \g_2^2} + \ddot \g_2^2 \right)\g_1\, dt 
	= \int_a^b  \eta \ddot \g_2^2\, \frac{\g_1}{\dot \g_1^2}\, dt \nonumber\\
		&  \geq \inf \left\{    \int_a^b  \dot \phi^2\, \frac{\g_1}{\dot \g_1^2}\, dt : \phi \in W^{1,2}_{\textup{loc}}((a,b);\R^2),\ \phi(a)=\dot \g_2(a),\ \right.\label{eq:minimize2}\\
		&\hspace{4cm} \left. \phi(b)=\dot \g_2(b),\ \dot \phi^2\, \frac{\g_1}{\dot \g_1^2}=0 \text{ on } \{\dot \g_1=0\}\right\}.\nonumber 
\end{align}
Let $\psi:= \g_1 / \dot \g_1^2$, $x_a:=\dot \g_2(a),\ x_b:=\dot \g_2(b)$. By the same computations as in \eqref{eq:defphi} and following lines, the minimum value is then given by
$$ \int_a^b \dot \phi^2(t)\, \psi(t)\, dt = \frac{(x_b-x_a)^2}{J},\qquad 
\text{with}\qquad
 J=\int_a^b \frac{1}{\psi(t)}\, dt = \int_a^b \frac{\dot \g_1^2(t)}{\g_1(t)}\, dt.$$
If $J=+\infty$, then \eqref{eq:oscillbd2} is trivially true. If $J<+\infty$, by inserting this minimum value in \eqref{eq:minimize2} and multiplying by $J$ we obtain
\begin{equation*}
 	\left(\frac{\ell}{2\pi}\int_a^b \k_1^2\, d\mu_\g\right)\left(\int_a^b \frac{\dot \g_1^2(t)}{\g_1(t)}\, dt\right) \geq |\dot \g_2(b)-\dot \g_2(a)|^2.
\end{equation*}
Using the simple inequality $ (x+y)^2 \geq 4xy$, we finally  obtain
\begin{equation*}
 	\ell \left(\int_a^b \k_1^2\, d\mu_\g + \int_a^b \frac{\dot \g_1^2(t)}{\g_1(t)}\, dt\right)^2 \geq 8\pi|\dot \g_2(b)-\dot \g_2(a)|^2.
\end{equation*}
Taking the square root, we obtain \eqref{eq:oscillbd2}.
\end{proof}

\begin{rem} In the case $\g_1\in W^{2,1}(a,b)$, $|\dot \g|=1$, $\dot\g_2\neq 0$ in $(a,b)$,
\begin{align*}
 	\frac1{4\pi} \int_a^b \k_1^2+ \k_2^2\, d\mu_\g &\geq \frac{1}{2\pi}\int_a^b |\k_1\k_2|\,d\mu_\g =  \int_a^b |\ddot \g||\dot \g_2|\, dt = \int_a^b |\dot \g_2| \sqrt{ \ddot \g_1^2 + \frac{\dot\g_1^2\ddot \g_1^2}{1-\dot \g_1^2} }\, dt\\
		& = \int_a^b |\dot \g_2| |\ddot \g_1|\sqrt{  \frac{ 1}{\dot \g_2^2} }\, dt = \int_a^b |\ddot \g_1|\, dt\geq |\dot \g_1(b)-\dot \g_1(a)|.
\end{align*}
We will make use of Lemma \ref{lemma:oscbd}, instead of the simpler estimate obtained in this Remark, because it allows for a rigorous treatment of the set $\{\dot \g_2=0\}$ and of the jump points for $\dot \g_1$.
\end{rem}

For $x\in \R$ we denote the integer part of $x$ by $\lfloor x \rfloor$.
\begin{lemma}
\label{lemma:finiteg10}
Let $\g$ be as in Definition \ref{def:gengen}. Then, there cannot be more than $\left\lfloor\frac{C}{8\pi}\right\rfloor$ intervals $(a_j,b_j)$ such that 
$$\dot \g_1(a_j)=\pm \ell(\g),\quad  \dot \g_1(b_j)=\mp \ell(\g).$$ 
As a consequence, 
$$
	\#\{\g_1=0\}\leq \left\lfloor\frac{C}{8\pi}\right\rfloor +1.
$$	
\end{lemma}
\begin{proof}
Let $(a,b)\subset [0,1]$ be an interval such that 
\begin{equation}
\label{eq:typint}
	\dot \g_1(a)= \ell(\g),\quad  \dot \g_1(b)=- \ell(\g)\quad \text{and}\quad \g_1>0\ \text{in }(a,b).
\end{equation}
By Corollary \ref{cor:reg} we have that $\g \in W^{2,2}_{\textup{loc}}\big((a,b);\R^2\big)$, $\dot \g \in C^0([a,b];\R^2)$. Since $|\dot \g_1(b)-\dot \g_1(a)|= 2\ell(\g)$, by Lemma \ref{lemma:oscbd}
$$ 
	\ell(\g) \int_a^b \k_1^2 +\k_2^2\, d\mu_\g \geq 4\pi |\dot \g_1(b)-\dot \g_1(a)|= 8\pi \ell(\g).
$$
Therefore, there cannot be more than $\left\lfloor\frac{C}{8\pi}\right\rfloor$ intervals satisfying \eqref{eq:typint}. In particular, by Lemma \ref{lemma:tg}, there cannot be more than $\left\lfloor\frac{C}{8\pi}\right\rfloor +1$ points where $\g_1=0$.
\end{proof}
The results obtained in this section imply that a generalized generator $\g$ with bounded Helfrich energy is either in (G1), if $\g_1>0$ and $\g$ is closed, or it can be decomposed into a finite number of curves in (G0). More precisely, we have the following result.

\begin{cor} 
\label{cor:split}
	If $\g$ as in Definition \ref{def:gengen} generates $\S$ and satisfies $\g_1(0)=\g_1(1)=0$, then there exist $k\in \N$ and curves $\eta_i$ generating $\S_i$ and satisfying  \eqref{eq:condreg}-\eqref{eq:condpos} for $i=1,\ldots,k$, such that 
\begin{equation}
\label{eq:eq}
	|\S| =\sum_{i=1}^k |\S_i|,\qquad \V(\S)=\sum_{i=1}^k\V(\S_i),\qquad \H(\S) =\sum_{i=1}^k \H(\S_i).
\end{equation}
Moreover, the generated surfaces $\S_i$ admit a $C^1$-regular parametrization.
\end{cor}	
\begin{proof}
By Lemma \ref{lemma:finiteg10}, $\#\{\g_1=0\}<+\infty.$ Therefore, there exist $0=t_0<t_1<\ldots<t_{k-1}<t_{k}=1$ such that $\{t_i\}_{i=0}^{k}=\{\g_1=0\}$. For $i=1,\ldots,k$, let $\ell_i:=t_i-t_{i-1}$ and define the curves
$$\eta_i:[0,1]\to \R^2,\qquad \eta_i(\t):=\g\left(\ell_i\t+ t_{i-1}\right).$$
It is immediate to check that $\eta_i$ satisfies \eqref{eq:condspeed}, \eqref{eq:condpos} and \eqref{eq:eq}, while \eqref{eq:condreg} is ensured by Corollary \ref{cor:reg}. The latter and Lemma \ref{lemma:tg} also imply that $\S_i$ has a  $C^1$-regular parametrization.
\end{proof}


\section{Existence of a minimizer}
\label{sec:existence}

\setcounter{equation}{0}
\setcounter{theo}{0}



\begin{defi}[\textit{Systems of generalized surfaces}] 
\label{def:systems}
We say that a system of surfaces $S$ belongs to the class $\SS$ of systems of generalized surfaces if 
\begin{itemize}
	\item there exists $m\in \N$ such that $S=(\S_1,\ldots,\S_m)$, 
	\item for each $i=1,\ldots,m$ there is a curve $\g_i$ as in Definition \ref{def:gengen} which generates $\S_i$,
	\item $(\g_i)\cap (\g_j)=\emptyset$ for all $i\neq j$.
\end{itemize}	  

\end{defi}
Note that the number of components may depend on the choice of parametrization. For example, a system of two spheres ${\mathbb S}^2$ touching in one point on the $z$-axis  can be parametrized by one generator $\g$ of length $|\dot \g|=2\pi$ or by two generators $\g_a,\g_b$ of length $|\dot \g_a|=|\dot \g_b|=\pi$.  In Section \ref{ssec:surrev} we proved that generalized generators $\g$ are piecewise-$C^1$ and are not differentiable only in a finite number of points on the $z$-axis. Therefore, the number of components we are interested in is the (minimum) number of $C^1$ components, that is
$$ \# S:= \sum_{i=1}^m (\# \{\g_{1,i}(t)=0,t\neq 0,t\neq 1\}+1 ).$$
We recall that area, volume, and Helfrich energy of a system $S=(\S_1,\ldots,\S_m)$ are simply defined as the sum of the correspondent \textit{generalized} quantities over all the components of the system, i.e.,
\begin{align} 
	|S|&=\sum_{i=1}^m |\S_i|= \sum_{i=1}^m 2\pi \int_0^1 \g_{1,i}(t)|\dot \g_i(t)|\, dt, & &	\F(S)=\sum_{i=1}^m \H(\S_i),\label{def:areasys}\\
	\V(S)&=\sum_{i=1}^m \V(\S_i)=\sum_{i=1}^m \pi \int_0^1 (\g_{1,i}(t))^2\dot \g_{2,i}(t)\, dt. & & \nonumber
\end{align}

\subsection{Convergence of measure-function couples}
\label{ssec:mt}
We turn now to the suitable notion of convergence for such systems. We recall that a sequence of Radon measures $\mu^n$ is said to converge weakly-$*$ to $\mu\in RM(\R)$ if
$$ \lim_{n\to \infty} \int_\R  \phi(t)\, d\mu^n(t) \to  \int_\R  \phi(t)\, d\mu(t)$$
for every $\phi \in C^0_c(\R)$.
We define the space of $p$-summable functions with respect to a positive Radon measure $\mu$ as
$$ L^p(\mu;\R^2):=\left\{ f:\R\to \R^2\ \mu\text{-measurable, such that }\int_\R |f(x)|^p\, d\mu(x)<+\infty \right\}.$$
 
\begin{defi}[\textit{Convergence of measure-function couples}]
\label{def:weakags}
Following \cite[Definition 5.4.3]{AGS}, given a sequence of measures $\mu^n\in RM(\R)$ converging weakly-$*$ to $\mu$, we say that a sequence of (vector) functions $f^n\in L^1(\mu^n;\R^2)$ converges weakly to a function $f\in L^1(\mu;\R^2)$, and we write $f^n \weakto f$ in $L^1(\mu^n;\R^2)$, provided
\begin{equation}
\label{eq:weakags}
	\lim_{n\to \infty} \int_\R f^n (t)\cdot \phi(t)\, d\mu^n(t) \to  \int_\R f (t)\cdot \phi(t)\, d\mu(t)
\end{equation}
for every $\phi \in C^\infty_c(\R;\R^2)$. For $p>1$, we say that a sequence of (vector) functions $f^n\in L^p(\mu^n;\R^2)$ converges \emph{weakly} to a function $f\in L^p(\mu;\R^2)$, and we write $f^n \weakto f$ in $L^p(\mu^n;\R^2)$, provided
\begin{equation}
\label{eq:weakagsp}
	\sup_{n\in \N} \int_\R |f^n (t)|^p\, d\mu^n(t) <+\infty\quad \text{and}\quad f^n \weakto f \text{ in }L^1(\mu^n;\R^2).
\end{equation}
For $p>1$, we say that a sequence of (vector) functions $f^n\in L^p(\mu^n;\R^2)$ converges \emph{strongly} to a function $f\in L^p(\mu;\R^2)$, and we write $f^n \to f$ in $L^p(\mu^n;\R^2)$, if \eqref{eq:weakagsp} holds and
$$ \limsup_{n \to \infty} {\|f^n\|}_{L^p(\mu^n;\R^2)}\leq {\|f\|}_{L^p(\mu;\R^2)}.$$
\end{defi}
\begin{lemma}[Weak-strong convergence in $L^p(\mu;\R^d)$ {\cite[Proposition 3.2]{Moser01}}]
\label{lemma:moser}
Let $p,q\in (1,\infty)$ such that $1/p+1/q=1$. Suppose that $\mu^n$ and $\mu$ are Radon measures on $\R$ such that $\mu^n\weaksto \mu$ and that $f^n\in L^p(\mu^n;\R^d)$, $f\in L^p(\mu;\R^d)$, $g^n\in L^q(\mu^n;\R^d)$, $g\in L^q(\mu;\R^d)$ be such that
$$
	f^n \weakto f\quad \text{weakly in }L^p(\mu^n;\R^d),\qquad 	g^n \to g\quad \text{strongly in }L^q(\mu^n;\R^d).
$$
Then
$$
	f^n g^n \weakto fg \quad \text{weakly in }L^1(\mu^n;\R^d).
$$
\end{lemma}
\begin{theo}[Lower-semicontinuity {\cite[Theorem 5.4.4 - (ii)]{AGS}}]
\label{th:ags}
Let $p>1$,  let $\mu^n$ and $\mu$ be Radon measures on $\R$ such that $\mu^n\weaksto \mu$, and let $f^n\in L^p(\mu^n;\R^2)$ be a sequence converging \emph{weakly} to a function $f\in L^p(\mu;\R^2)$ in the sense of Definition \ref{def:weakags}. Then
$$
	\liminf_{n \to \infty} \int_\R g(f^n(t))\, d\mu^n(t) \geq \int_\R g(f(t))\, d\mu(t),
$$
for every convex and lower-semicontinuous function $g:\R\to (-\infty,+\infty].$
\end{theo}
\begin{defi}[\textit{Convergence of systems of surfaces}] 
\label{def:systconv}
We say that a sequence of systems of surfaces $S^n\in \SS$ as in Definition \ref{def:systems} converges to a finite system $S\in \SS$ if $S$ admits a parametrization $\g_1,\ldots,\g_w$ and there exist $\bar n \in \N$ and $\bar m\geq w$ such that $\# S^n =\bar m$ for all $n\geq \bar n$ and such that 
\begin{itemize}
	\item[(i)] for each $i=1,\ldots,w$ the generating curves $\g^n_i$ converge in the following sense
		\begin{align}
			\g^n_i &\to \g_i \quad \text{uniformly in }C^0([0,1];\R^2), 									\label{eq:convg}\\
			\dot \g^n_i &\to  \dot \g_i\quad \text{strongly in } L^2\val, 							\label{eq:convgd}\\
			\ddot \g^n_i &\weakto \ddot \g_i \quad \text{weakly in $L^2(\mu_{\g^n};\R^2)$, in the sense of \eqref{eq:weakagsp}};	\label{eq:convgdd}
		\end{align}
	\item[(ii)] for each $i=w$+$1,\dots, \bar m$ 
		$$ \g^n_i \to 0\quad \text{strongly in }W^{1,2}(0,1;\R^2).$$
\end{itemize}		
\end{defi}
Note that \eqref{eq:convg} and \eqref{eq:convgd} imply that $\mu_{\g^n}\weaksto \mu_\g$, so that \eqref{eq:convgdd} is well-defined.

\subsection{Compatibility of constraints}
\label{ssec:constr}
With this definition of convergence, the passage to the limit of the area and volume constraints is straightforward. Let $S^n \in \mathcal A(A,V)$ be a sequence of systems converging to a system $S=(\S_1,\ldots,\S_w)\in \SS$ in the sense of Definition \ref{def:systconv}, then
\begin{align}
	&\lim_{n \to\infty} |S^n|=|S|,\label{eq:areasconv}\\ 
	&\lim_{n \to\infty} \V(S^n)=\V(S). \label{eq:vsconv}
\end{align}	
Indeed, by \eqref{eq:convg} and \eqref{eq:convgd}, using the fact that the for the last $(\bar m-w)$ components $\dot \g^n_i\to 0$ strongly in $L^2(0,1)$,
\begin{equation*}
	A=\lim_{n\to \infty}\sum_{i=1}^{\bar m}|\S_i^n|= \lim_{n \to \infty}\sum_{i=1}^{\bar m} 2\pi \int_0^1 \g_{i,1}^n |\dot \g_i^n|\, dt= \sum_{i=1}^w2\pi \int_0^1 \g_{i,1} |\dot \g_i|\, dt = \sum_{i=1}^w|\S_i|= |S|.
\end{equation*}
and
\begin{align*}
	V=\lim_{n\to \infty}\sum_{i=1}^{\bar m}\V(\S_i^n) &= \lim_{n \to \infty}\sum_{i=1}^{\bar m} \pi \int_0^1 (\g_{i,1}^n)^2 \dot \g_{i,2}^n\, dt\\
		&= \sum_{i=1}^w\pi \int_0^1 \g_{i,1}^2 \dot \g_{i,2}\, dt = \sum_{i=1}^w\V(\S_i)= \V(S). 
\end{align*}


\subsection{Lower-semicontinuity}
\label{subsec:lsc}
\begin{prop}
\label{prop:lscS} 
Let $S,S^n\in \SS $ be a family of systems of surfaces such that $ S^n \to S$  $\text{in the sense of Definition \ref{def:systconv}}. $
Then
$$ \liminf_{n \to \infty} \F(S^n)\geq \F(S).$$
\end{prop}
\begin{proof}
Let $S=(\S_1,\ldots,\S_w)$. First of all, we notice that, according to Definition \ref{def:systconv}, since $\F\geq 0$, for all $n$ big enough
$$  \F(S^n) = \sum_{i=1}^w \F(\S_i^n) + \sum_{i=w+1}^{\bar m}\F(\S_i^n) \geq \sum_{i=1}^w \F(\S_i^n).$$
Therefore, it is enough to show that
$$ \liminf_{n \to \infty} \F(\tilde S^n)\geq \F(S)$$
for $\tilde S^n=(\S^n_1,\ldots,\S^n_w)$. This is reasonable, since by Definition \ref{def:systconv}-(ii) the $i$-th component, for $i>w$, vanishes as measure.

Let then $i\in \{1,\ldots,w\}$ be fixed, we write $\g^n=\g^n_i$ for the generating curve and $\S^n=\S^n_i$ for the generated revolution surface. Define the family of measures 
$$
	\mu_{\g^n}:=2\pi\, \g_1^n|\dot \g^n|\, \L^1|_{[0,1]} \in RM(\R).
$$
By \eqref{eq:convg} and \eqref{eq:convgd}
\begin{equation}
\label{eq:muconv}
	\mu_{\g^n} \weaksto \mu :=2\pi\, \g_1|\dot \g| \L^1|_{[0,1]}\qquad \text{in } RM(\R).
\end{equation}
Recall that 
\begin{align*}
	\k_1&=\frac{\ddot \g_2 \dot \g_1 -\ddot \g_1 \dot \g_2}{|\dot \g|^3}, &&	|\k_1|=\frac{|\ddot \g|\ }{|\dot \g|^2},		&&
	\k_2=\frac{\dot \g_2}{\g_1 |\dot \g| }
\end{align*}	
($\k_1^n$ and $\k_2^n$ are defined analogously) and
\begin{equation}
\label{eq:4integrals}
	\H(\S^n) = \int_\R \left\{\frac{\kh}{2}(\k_1^n +\k_2^n-H_0)^2 +\kg k_1^n k_2^n\right\}\, d\mu_{\g^n}.
\end{equation}
Let $\S^{n_k}$ be a subsequence realizing the liminf of $\H$, i.e., such that
$$ \lim_{k\to \infty} \H(\S^{n_k})=\liminf_{n\to \infty}\H(\S^n).$$
For sake of notation, we drop the index $k$ from the subsequence in the rest of the proof. All we need to prove is that 
\begin{align}
	&\k_1^n \weakto \k_1\qquad \text{weakly in }L^2(\mu_{\g^n};\R), \label{eq:convk1}\\
	&\k_2^n \to \k_2\qquad \text{strongly in }L^2(\mu_{\g^n};\R), \label{eq:convk2}
\end{align}
in the sense of Definition \ref{def:weakags}. Indeed, if \eqref{eq:convk1} and \eqref{eq:convk2} hold,  by Theorem \ref{th:ags} and \eqref{eq:muconv}
\begin{equation*}
	\liminf_{n\to \infty} \int_\R (\k_1^n +\k_2^n -H_0)^2\,  d\mu_{\g^n}   \geq \int_\R (\k_1 +\k_2 -H_0)^2\,  d\mu 
\end{equation*}
and by Lemma \ref{lemma:moser}
\begin{equation*}
	\lim_{n\to \infty} \int_\R \k_1^n \k_2^n \,  d\mu_{\g^n}   = \int_\R \k_1 \k_2 \,  d\mu. 
\end{equation*}
Therefore,  for every component $\S^n_i$ 
\begin{equation}
\label{eq:liminfs}
	\liminf_{n \to \infty} \H(\S^n_i) \geq \H(\S_i),
\end{equation}
and to conclude the proof, it is sufficient to notice that
$$
	\liminf_{n \to \infty} \F(\tilde S^n) = \liminf_{n \to \infty} \sum_{i=1}^w \H(\S^n_i) \geq \sum_{i=1}^w \liminf_{n \to \infty} \H(\S^n_i) \stackrel{\eqref{eq:liminfs}}{\geq} \sum_{i=1}^w \H(\S_i) =\F(S).
$$	

Convergences  \eqref{eq:convk1}-\eqref{eq:convk2} are addressed in the next Lemma, which concludes the proof of Proposition \ref{prop:lscS}.
\end{proof}

\begin{lemma}
\label{lemma:k1ags}
Let $\g^n:[0,1]\to\R^2$ be a sequence of generating curves for admissible surfaces $\S^n$, and assume that $\g^n\to\g$ as in \eqref{eq:convg}-\eqref{eq:convgdd}. Then 
\begin{align*} 
	\k_1^n&=\frac{\ddot \g_2^n\dot \g_1^n-\ddot \g_1^n\dot \g_2^n}{|\dot \g^n|^3}\hspace{1cm} \text{converges weakly to }\k_1 = \frac{\ddot \g_2\dot \g_1-\ddot \g_1\dot \g_2}{|\dot \g|^3},\\
	\k_2^n&=\frac{\dot \g_2^n}{\g_1^n |\dot \g^n| }\hspace{2.1cm} \text{converges strongly to } \k_2=\frac{\dot \g_2}{\g_1 |\dot \g| },
\end{align*}
in the sense of Definition \ref{eq:weakags}, where $\k_1$ and $\k_2$ are defined $\mu$-a.e.
\end{lemma}
\begin{proof}
Note that, in case (i), we can assume that
\begin{align}
	|\g^n|	&\leq M<\infty,		\label{eq:case1m}\\
	|\dot \g^n|&\geq L >0.	\label{eq:case1l}
\end{align}	
Denote by $R\psi$ the $\pi/2$ rotation of the function $\psi=(\psi_1,\psi_2)$, i.e., $R\psi=(-\psi_2,\psi_1)$. Then, for all $\phi\in C^1_c(\R)$
\begin{equation*}
	\int_\R \k_1^n \phi\, d\mu_{\g^n}= \int_\R \left(\ddot \g^n\cdot \frac{R\dot \g^n}{|\dot \g^n|^3}\right)\phi\, d\mu_{\g^n}.
\end{equation*}
By \eqref{eq:convg} and \eqref{eq:convgd}
\begin{equation*}
	\lim_{n\to \infty} \int_\R \left(\frac{R\dot \g^n}{|\dot \g^n|^3}\right)^2 d\mu_{\g^n} = \lim_{n\to \infty} 2\pi \int_\R \frac{1}{|\dot \g^n|^3} \g^n_1\,dt =2\pi \int_\R \frac{1}{|\dot \g|^3} \g_1\,dt = \int_\R \left(\frac{R\dot \g}{|\dot \g|^3}\right)^2 d\mu_{\g},
\end{equation*}
therefore, $R\dot \g^n /|\dot \g^n|^3$ converges strongly in $L^2(\mu_{\g^n};\R^2)$ to $R\dot \g /|\dot \g|^3$, according to Definition \ref{def:weakags}. By \eqref{eq:convgdd} and Lemma \ref{lemma:moser} 
\begin{align*}
	\lim_{n\to \infty} \int_\R \k_1^n \phi\, d\mu_{\g^n} &= \lim_{n\to \infty}\int_\R \left(\ddot \g^n\cdot \frac{R\dot \g^n}{|\dot \g^n|^3}\right)\phi\, d\mu_{\g^n}\\
			& =\int_\R \left(\ddot \g\cdot \frac{R\dot \g}{|\dot \g|^3}\right)\phi\, d\mu_{\g} = \int_\R \k_1 \phi\, d\mu_{\g}.
\end{align*}	
Regarding $\k_2$,  let $\phi\in C^1_c(\R)$, then
\begin{align*}
	\int_\R \k_2^n(t)\phi(t)\, d\mu^n(t)= 2\pi \int_0^1 \frac{\dot \g_2^n(t)}{\g_1^n|\dot \g^n|}\phi(t)\, \g_1^n|\dot \g^n|dt = 2\pi \int_0^1 \dot \g_2^n(t) \phi(t)\, dt,
\end{align*}
so that, by \eqref{eq:convgd},
\begin{equation}
\label{eq:k2ags}
	\lim_{n\to \infty} \int_\R \k_2^n(t)\phi(t)\, d\mu^n=2\pi \int_0^1 \dot \g_2(t) \phi(t)\, dt= \int_\R \k_2(t)\phi(t)\, d\mu(t).
\end{equation}
By \eqref{eq:convg} and \eqref{eq:convgd} we also have that
\begin{align*}
	\lim_{n \to \infty} {\|\k_2^n\|^2}_{L^2(\mu_{\g^n};\R^2)} &= \lim_{n\to \infty} \int_\R (\k_2^n(t))^2\, d\mu^n = \lim_{n\to \infty} 2\pi \int_0^1 \frac{(\dot \g_2^n)^2}{\g_1^n}\, dt\\
	 &= 2\pi \int_0^1 \frac{(\dot \g_2)^2}{\g_1}\, dt = \int_\R k_2^2\, d\mu_\g = {\|k_2\|^2}_{L^2(\mu_{\g};\R^2)}. 
\end{align*}
This concludes the proof of Lemma \ref{lemma:k1ags} which completes the proof lower-semicontinuity statement.
\end{proof}
As a direct consequence of Lemma \ref{lemma:k1ags} and Theorem \ref{th:ags} it holds
\begin{cor}
\label{cor:separatelsc}
Under the assumptions of Proposition \ref{prop:lscS}, by Lemma \ref{lemma:k1ags} and Theorem \ref{th:ags}
\begin{equation}
\label{eq:seplsc}
	\liminf_{n\to \infty} \int_\R(\k_1^n)^2+(\k_2^n)^2\,d\mu_{\g^n} \geq \int_\R (\k_1)^2+(\k_2)^2\,d\mu_{\g}.
\end{equation}
\end{cor}


\subsection{Compactness} 
\label{ssec:compact}

Let $\{S^n\}$ be a family of finite systems generated by curves in (G0) or (G1). Throughout this step, we assume that there exist constants $\Lambda,A$ such that for all $n\in \N$
\begin{equation}
\label{assumptions:bounds} 
	|S^n|\leq A,\qquad \F(S^n)\leq \Lambda.
\end{equation}	
Moreover, owing to Lemma \ref{lemma:lengthbound} and \eqref{eq:totalbound}, the diameter of every component of a system satisfying \eqref{assumptions:bounds} is bounded by a constant which depends only on the data of the problem. Consequently, since we are interested only in the shape of the components, and not in their relative position in the space, it is not restrictive to assume that there exists $R>0$ such that for all $n\in \N$
\begin{equation}
\label{assumptions:bdm}
	\S^n \subset B_R.
\end{equation}
For short, we say that a finite system $S$ generated by curves in (G0) and (G1) satisfying \eqref{assumptions:bounds} and \eqref{assumptions:bdm} is an \textit{admissible system}. Note that in \eqref{assumptions:bounds} we are not fixing the total area, we are requiring only an upper bound, thus allowing for more flexibility in the proofs. Note also that owing to the isoperimetric inequality no separate bound on the volume is necessary, and owing to the restrictive geometry of revolution surfaces the only possible genus is 0 or 1.


\begin{lemma} 
\label{lemma:fincard}
Let $S=(\S_1,\ldots,\S_m)$ be a finite system generated by curves in \GO or \GI and satisfying \eqref{assumptions:bounds}, then there exists a constant $C=C(H_0,\kh,\kg,A,\Lambda)$ such that $\# S=m\leq C$.
\end{lemma}

\begin{proof}
If $\g$ is in the class (G0), by Lemma \ref{lemma:tg} we have 
\begin{equation}
\label{eq:bdg1}
	|\dot \g_1(1)-\dot\g_1(0)|=2\ell(\g).
\end{equation}
If $\g$ is in (G1), then  $\g\in C^1([0,1])$ by Corollary \ref{cor:reg}, and therefore there must be two points $s,t \in [0,1]$ such that 
\begin{equation}
\label{eq:bdg2}
	|\dot \g_1(s)-\dot\g_1(t)|=2\ell(\g)
\end{equation}
(e.g, choose $s\in \text{argmin}(\g_2)$ and $t\in \text{argmax}(\g_2)$: since $\g$ is closed and differentiable, the tangents in these points are horizontal and have opposite orientation). 
By  \eqref{eq:bdg1}, \eqref{eq:bdg2} and Lemma \ref{lemma:oscbd}, for every $i=1,\ldots,m$
$$ 
	8\pi \ell(\g_i)\leq 4\pi \max_{s,t\in[0,1]}|\dot \g_1(s)-\dot \g_1(t)|	\leq \ell(\g_i)\int_0^1  \k_{1,i}^2+ \k_{2,i}^2\, d\mu_{\g_i},
$$
$$ 
	8\pi \leq  \int_0^1  \k_{1,i}^2+ \k_{2,i}^2\, d\mu_{\g_i}. 
$$ 	
By \eqref{eq:totalbound} there is a constant $C_1>0$, depending only on the data, such that summing over $i$
$$
	8\pi m = \sum_{i=1}^m 8 \pi \leq \sum_{i=1}^m\int_0^1  \k_{1,i}^2+ \k_{2,i}^2\, d\mu_{\g_i}\leq C_1 (\F(S)+|S|).
$$
We conclude that $m\leq C_1(\Lambda +A)/8\pi:=C$.
\end{proof}

Since every admissible sequence $S^n$ satisfies $1\leq \# S^n\leq C$, we can now study the compactness property for each component separately. From here onwards we omit the component index $i$.


\begin{prop}[Compactness]
\label{prop:compsurf}
Let $\g^n:[0,1]\to\R^2$ be a sequence of generating curves for admissible surfaces $\S^n$. Then, either
\begin{itemize} 
	\item[(i)] there exists a subsequence $\g^{n_k}$ and a generalized generator $\g$ as in Definition \ref{def:gengen} such that $\g^{n_k}$ converges to $\g$ in the sense of \eqref{eq:convg}-\eqref{eq:convgdd}, 

	or
	\item[(ii)] there exists a point $(0,z)\in \R^2$ such that $\g^{n} \to (0,z)$ strongly in $W^{1,2}((0,1);\R^2)$.
\end{itemize}
\end{prop}
\begin{proof}
Let $\S^n$ be a sequence of admissible surfaces generated by $\g^n:[0,1]\to\R^2$, parametrized so that $|\dot \g(t)|\equiv \ell(\g^n)$. It holds 
\begin{equation}
\label{eq:tb3} 
	\sup_{n\in \N} \int_0^1 (\k_1^n)^2+(\k_2^n)^2\,d\mu_{\g^n}\stackrel{\eqref{eq:totalbound}}{\leq} C(\H(\S^n)+|S^n|)\stackrel{\eqref{assumptions:bounds}}{\leq} C(\Lambda +A).
\end{equation}

\paragraph{Step I.} Proof of case (ii) and convergence \eqref{eq:convg} in case (i). By \eqref{assumptions:bdm}, Lemma \ref{lemma:lengthbound} and \eqref{eq:tb3}, there exists a constant $C>0$ such that 
\begin{equation}
\label{eq:prop1:00} 
	\max_{t\in [0,1]}|\g^n(t)| + 	|\dot \g^n(t)|\leq \max_{t\in [0,1]}|\g^n(t)| + |\ell(\g^n)|\leq C
\end{equation}
for every $n\in \N$.	 Therefore, by Arzel\`a-Ascoli Theorem and weak-$*$ compactness in $L^\infty$ there exists a subsequence, which we do not relabel, and a there exists a continuous limit curve $\g:[0,1]\to \R^2$ such that
\begin{align}
	\g^n &\to \g \quad \text{ uniformly in }C^0([0,1];\R^2),\label{eq:prop1:0}\\
	\dot \g^n &\weaksto  \dot \g\quad \text{ weakly-$*$ in } L^\infty((0,1);\R^2)\label{eq:prop1:1}.
\end{align}
Moreover, up to extracting a further subsequence, since $\{|\dot \g^n|\}$ is just a bounded sequence of real numbers, we can assume that  $ 	|\dot \g^n| \to L$. The point now is to prove that $|\dot \g|=L$.
	
If $L=0$, by lower-semicontinuity of the $L^2$-norm  with respect to the weak topology
\begin{equation}
\label{eq:tozero}
	0 = \lim_{n\to\infty}\int_0^1|\dot \g^n|^2\, dt \geq \int_0^1|\dot \g|^2\, dt \geq 0,
\end{equation}
i.e., $\dot \g^n \to 0$ strongly in $L^2((0,1);\R^2)$ and $\g^n$ strongly converges in $W^{1,2}((0,1);\R^2)$ to a constant $\g$. If, by absurd $\g_1>0$, then there exists $\e>0$ such that $\g_1^n(t)>\e$ for all $n$ big enough and thus (as in \cite[Lemma 3.1]{BellDMasoPaolini93}) 
$$ (2\pi)^2\leq \ell(\g^n)\int_0^{\ell(\g^n)} |\ddot \g^n|^2\, dt \leq \frac{\ell(\g^n)}{\e}\int_0^{\ell(\g^n)} |\ddot \g^n|^2\, \g_1^n dt\leq \frac{\ell(\g^n)}{\e}\int_{\S^n} (k_1^n)^2 dA, $$
and therefore, $\ell(\g^n)>\frac{4\pi^2\e}{\Lambda}$, which contradicts \eqref{eq:tozero}. This covers case $(ii)$ of Proposition \ref{prop:compsurf}. We also note that
\begin{equation*}
	\lim_{n \to \infty}|\S^n|= 0\qquad \text{if and only if}\qquad  \lim_{n \to \infty}|\dot \g^n| =0.
\end{equation*}
Indeed, the ``only if" part follows directly from Lemma \ref{lemma:lengthbound} and \eqref{eq:tb3}, while the ``if" part is a consequence of  \eqref{eq:defarea}, \eqref{eq:prop1:0} and \eqref{eq:tozero}.  In the following steps, in order to prove case $(i)$, we assume
\begin{equation}
\label{eq:prop1:2}
	|\dot \g^n| \to L>0.
\end{equation}

\paragraph{Step II.} Convergence \eqref{eq:convgd} in case (i). We have to show that 
$$ \dot \g^n \to  \dot \g\quad \text{ strongly in } L^2\val. $$
Recall that the Total Variation of a function $f:[0,1]\to \R$ is defined as
$$ V_{[0,1]}(f):=\sup\left\{ \sum_{j=0}^{N-1}|f(t_{i+1})-f(t_i)| : \{ 0=t_0, \ldots, t_N=1\} \text{ is a partition of }[0,1] \right\}. $$
By Lemma \ref{lemma:oscbd}, for any interval $(a,b)\subseteq (0,1)$
$$  	4\pi|\dot \g_1^n(b)-\dot \g_1^n(a)|\leq |\dot \g^n|\int_a^b (\k_1^n)^2+ (\k_2^n)^2\, d\mu_{\g^n}. $$
Therefore, by \eqref{eq:prop1:2} and \eqref{eq:tb3}, there exists a constant $C>0$ such that
$$ V_{[0,1]}(\dot \g_1^n)\leq \frac{|\dot \g^n|}{4\pi}\int_0^1 (\k_1^n)^2+ (\k_2^n)^2\, d\mu_{\g^n}\leq C.$$
By standard compactness in Bounded Variation spaces (see, e.g., \cite{EvansGariepy92}), there exists a subsequence $\g_1^{n_k}$ and a limit function $\s\in BV((0,1);\R^2)$ such that
$$ \dot \g_1^{n_k} \to  \s\quad \text{ strongly in } L^1(0,1). $$
Since $\dot \g_1$ is also bounded in $L^\infty(0,1)$
$$ \dot \g_1^{n_k} \to   \s\quad \text{ strongly in } L^p(0,1)\quad \forall\, p\in [1,+\infty), $$
and by \eqref{eq:prop1:1} we can identify the weak limit $\dot \g_1$ with the strong limit $\s$. We have thus proved
$$ \dot \g_1^n \to  \dot \g_1 \quad \text{ strongly in } L^p(0,1)\quad \forall\, p\in [1,+\infty). $$
Since, by hypothesis, $\g^n_1(t)> 0$ for all $t\in (0,1)$,  $\g_1\geq 0$ for all $t\in [0,1]$. Let
$$
	E:=\{t \in [0,1]:\g_1(t)=0\}.
$$
Note that $E$ is a compact subset of $[0,1]$, so, in particular, it is Lebesgue-measurable. Assume, by contradiction, that $|E|>0.$ Owing to the bound on the energy, there is a constant $C>0$ such that
$$
	C \geq \int_{\S^n} (\k_2^n)^2\, dA \geq  2\pi \int_0^1 \frac{(\dot \g_2^n)^2}{\g_1^n|\dot \g^n|}dt  \geq  2\pi \int_E \frac{(\dot \g_2^n)^2}{\g_1^n|\dot \g^n|}dt  \geq  \frac{2\pi}{|\dot \g^n|} \int_E \frac{(\dot \g_2^n)^2}{\e}dt.
$$
By \eqref{eq:prop1:2} $|\dot\g^n|$ is bounded from above and we can find a new constant $C$ such that
$$ 
	\forall\, \e>0\ \exists\, \bar n\in \N:\quad  {\|\dot \g^n_2\|}^2_{L^2(E)}\leq \e C \qquad \forall\, n>\bar n,
$$
which implies that $\dot \g_2^n\to 0$, strongly in $L^2(E)$. Since $(\dot \g_1^n)^2=\ell(\g^n)^2-(\dot \g_2^n)^2$, by \eqref{eq:prop1:2} $|\dot \g_1^n|\to  L$ strongly in $L^2(E)$. Since strong convergence implies pointwise convergence (up to extracting a subsequence), we obtain that $\dot \g_1(t) \in \{-L,L\}$ for a.e. $t\in E$. On the other hand, since every $t\in E$ is a point of minimum for $ \g_1$, it must be $\dot \g_1(t)=0$ for a.e. $t\in E$. We have thus obtained a contradiction and we conclude that $|E|=0$.

For $h \in \N$ define
$$ A_h:=\left\{ t \in [0,1]:\g_1(t)>\frac 2h \right\}.$$ 
Since $A_h \subset A_{h+1}$ and $[0,1]\backslash E \subseteq \cup_{h \in \N}A_h$, we have that 
\begin{equation}
\label{eq:measah}
	\lim_{h \to \infty}|[0,1]\backslash A_h|=0.
\end{equation}
Our aim is now to show that 
\begin{equation}
\label{eq:convah}
	\dot \g_2^n \to  \dot \g_2 \quad \text{ strongly in } L^1(A_h).
\end{equation}
In fact, since
$$ \int_0^1|\dot \g_2^n -\dot \g_2|dt \leq \int_{A_h}|\dot \g_2^n -\dot \g_2|dt +(|\dot \g^n| +L)|[0,1]\backslash A_h|,$$
then \eqref{eq:measah} and \eqref{eq:convah} imply that $ \dot \g_1^n \to  \dot \g_1 \text{ strongly in } L^1(0,1)$ and thus
$$ \dot \g_1^n \to  \dot \g_1 \quad \text{ strongly in } L^p(0,1),\quad \forall\, p\in [1,+\infty). $$
In order to prove \eqref{eq:convah}, note that since $\g_1^n$ converges uniformly to $\g_1$, for every $h\in \N$ there is $\bar n \in \N$ such that 
\begin{equation}
\label{eq:g1h}
	\g_1^n(t)>\frac 1h\quad \forall\,t\in A_h,\quad \forall\, n> \bar n. 
\end{equation}
Let now $h$ be fixed, $\bar n$ given as above and let $(a,b)\subset A_h$. By \eqref{eq:oscillbd2}, for all $n>\bar n$
\begin{align*}  	
	2\sqrt{2\pi}|\dot \g_2^n(b)-\dot \g_2^n(a)| & \leq |\dot \g^n|^{1/2}\left(\int_a^b (\k_1^n)^2 d\mu_{\g^n}+ \int_a^b\frac{(\dot \g_1^n)^2}{\g_1^n}\, dt\right)\\
			&\stackrel{\eqref{eq:g1h}}{\leq} |\dot \g^n|^{1/2}\left(\int_a^b (\k_1^n)^2 d\mu_{\g^n}+ h \int_a^b (\dot \g_1^n)^2\, dt\right).
\end{align*}	
As we did for $\dot \g_1^n$, we can then control $ V_{\overline A_h}(\dot \g_2^n)\leq C(1+h),$ and conclude \eqref{eq:convah}. Strong convergence of $\dot \g_1^n$ and $\dot \g_2^n$
, in particular, implies
\begin{equation}
\label{eq:strongconvdot}
	\lim_{n\to \infty}|\dot \g^n| = |\dot \g|.
\end{equation}

\paragraph{Step III.} Convergence \eqref{eq:convgdd} in case (i). The proof of convergence \eqref{eq:convgdd} follows from the strong-$L^2$ convergence of $\dot \g^n$.  Integrating by parts on $(0,1)$, for all $\phi \in C^1_c(\R)$
\begin{align*}
	\int_\R \ddot \g^n_1 (t) \phi(t)\, d\mu_{\g^n}(t) &= \frac{1}{|\dot \g^n|^3}\int_0^1 \frac{d}{dt}\left\{\dot \g^n_1 (t)\right\} \phi(t)\, \g^n_1(t)\, dt  \\
			&= -\frac{1}{|\dot \g^n|^3}\left\{ \int_0^1 (\dot \g^n_1 (t))^2 \phi(t)\, dt+ \int_0^1 \dot \g^n_1 (t)\g^n_1 (t) \dot \phi(t)\, dt\right\}.
\end{align*}
By \eqref{eq:prop1:0} and \eqref{eq:strongconvdot} we can pass to the limit as $n \to \infty$ and integrate back by parts
\begin{align}	
	\lim_{n\to \infty} \int_\R \ddot \g^n_1 (t) \phi(t)\, d\mu_{\g^n}(t) &= -\frac{1}{|\dot \g|^3}\left\{ \int_0^1 (\dot \g_1 (t))^2 \phi(t)\, dt+ \int_0^1 \dot \g_1 (t)\g_1 (t) \dot \phi(t)\, dt\right\}  \nonumber\\
	&=  \frac{1}{|\dot \g|^3}\int_0^1 \ddot \g_1 (t) \phi(t)\, \g_1(t)\, dt \nonumber\\
	&=  \int_\R \ddot \g_1 (t) \phi(t)\, d\mu_{\g}(t).\label{eq:gddot1}
\end{align}
In the same way, we get that
\begin{equation}
\label{eq:gddot2}
	\lim_{n\to \infty} \int_\R \ddot \g^n_2 (t) \phi(t)\, d\mu_{\g^n}(t)= \int_\R \ddot \g_2 (t) \phi(t)\, d\mu_{\g}(t).
\end{equation}
Limits \eqref{eq:gddot1} and \eqref{eq:gddot2} prove convergence \eqref{eq:convgdd}. 

%
%
\end{proof}

We now have all the ingredients to prove Theorem \ref{th:main}.

\subsection{Proof of Theorem \ref{th:main}}
\label{ssec:proof}
Let the area and volume constraints $A,V$ be given, such that the isoperimetric inequality \eqref{eq:isop} is satisfied. Let the set $\A(A,V)$ and the functional $\F$ be given as in the statement of Theorem \ref{th:main}. Let $S^n=(\S_1^n,\ldots,\S_{m(n)}^n)\in \A(A,V)$, where $m(n)=\# S^n$, be a sequence of systems of revolution surfaces such that
\begin{equation}
\label{eq:minimizing}
	\lim_{n\to \infty}\F(S^n)=\inf_{S\in \A(A,V)}\F(S).
\end{equation}

Note that by \eqref{eq:totalbound}
$$\sum_{i =1}^{m(n)}\left(\frac{|\S_i^n|}{2} +\int_{\S_i^n} (\k_{1,i}^n)^2+(\k_{2,i}^n)^2\,dA\right) \leq C(A+\F(S^n)).$$ 
Since $S^n$ satisfies $\F(S^n)\leq \Lambda$, for a suitable $\Lambda>0$, by Lemma \ref{lemma:fincard} there is a constant $C>0$ such that $\# S^n<C$ for all $n\in \N$. We can thus extract a subsequence (not relabeled) and find $m\in \N$, such that  $\#S^n \equiv m$ for all $n$. By Corollary \ref{cor:realbd} the total length of the curves generating $S^n$ is uniformly bounded and therefore it is not restrictive to assume that there exists $R>0$ such that every $S^n$ is contained in the ball of radius $R$ centered at the origin. Thus $S^n$ is \textit{admissible} in the sense of Section \ref{ssec:compact}. We can therefore apply the compactness result, i.e. Proposition \ref{prop:compsurf}. Up to extracting $m$ subsequences, we find generalized generators $\g_1,\ldots,\g_j$ and real numbers $z_1,\ldots,z_{m-j}$, for $0<j\leq m$, such that, as $n\to \infty$  
$$ (\g_1^n,\ldots,\g_j^n)\to (\g_1,\ldots,\g_j),$$
in the sense of convergence \eqref{eq:convg}-\eqref{eq:convgdd} and
$$ (\g_{j+1}^n,\ldots,\g_m^n)\to ((0,z_1),\ldots,(0,z_{m-j}))$$
strongly in $W^{1,2}((0,1);\R^2)$. Denoting by $S$ the system of surfaces generated by $(\g_1,\ldots,\g_j)$, by the lower-semicontinuity Proposition \ref{prop:lscS}
$$ \liminf_{n \to \infty} \F(S^n)\geq \F(S),$$
so that, by \eqref{eq:minimizing}, $\F(S)= \inf \F$. To conclude, we have to show that $S\in \A(A,V)$, or, more precisely, that $S$ can be (re)parametrized  by curves in (G0) or (G1), which satisfy the area and volume constraints. By Section \ref{ssec:constr}, the constraints are continuous with respect to strong convergence in $W^{1,2}((0,1);\R^2)$, so that $|S|=A$ and $\V(S)=V.$ Let $\g$ be a generator of $S$, we distinguish two cases. If $\g_1(t)>0$ for all $t\in [0,1]$, then $\g$ is the limit of a sequence $\g^n\in$ (G1), which implies that $\g$ is closed;  by Lemma \ref{lemma:w2loc} $\g$ satisfies \eqref{eq:condreg1} and we infer that $\g\in$ (G1).  In the second case, $\{\g_1=0\}\neq \emptyset$. Since it's not restrictive to assume that $\g_1(0)=\g_1(1)=0$,  we can apply Corollary \ref{cor:split} and conclude that $\g$ can be reparametrized as a union of curves in (G0), which generate a family of surfaces with the same area, enclosed volume, and Helfrich energy as the one generated by $\g$. The proof of Theorem \ref{th:main} is thus complete.

As a final note, recalling Remark \ref{rem:index} on the index of a system of curves, if we restrict the minimization to the class of systems of disjoint curves such that $I(S,p)\in \{0,1\}$ for almost every $p\in \R^2$, by continuity of the index under uniform convergence, we obtain that the index takes values 0 or 1 also for the limit system. In particular, this implies that the minimizer is without self-crossings.

 \subsubsection*{Acknowledgements} We are deeply indebted to Eliot Fried, for introducing us to this subject and for many stimulating discussions. M.\,V. would also like to thank Luca Mugnai for all the technical and spiritual support.

\end{document}